

\ifx\shlhetal\undefinedcontrolsequence\let\shlhetal\relax\fi

\input amstex
\expandafter\ifx\csname mathdefs.tex\endcsname\relax
  \expandafter\gdef\csname mathdefs.tex\endcsname{}
\else \message{Hey!  Apparently you were trying to
  \string\input{mathdefs.tex} twice.   This does not make sense.} 
\errmessage{Please edit your file (probably \jobname.tex) and remove
any duplicate ``\string\input'' lines}\endinput\fi




\catcode`\X=12\catcode`\@=11

\def\n@wcount{\alloc@0\count\countdef\insc@unt}
\def\n@wwrite{\alloc@7\write\chardef\sixt@@n}
\def\n@wread{\alloc@6\read\chardef\sixt@@n}
\def\r@s@t{\relax}\def\v@idline{\par}\def\@mputate#1/{#1}
\def\l@c@l#1X{\firstpart.#1}\def\gl@b@l#1X{#1}\def\t@d@l#1X{{}}

\def\crossrefs#1{\ifx\all#1\let\tr@ce=\all\else\def\tr@ce{#1,}\fi
   \n@wwrite\cit@tionsout\openout\cit@tionsout=\jobname.cit 
   \write\cit@tionsout{\tr@ce}\expandafter\setfl@gs\tr@ce,}
\def\setfl@gs#1,{\def\@{#1}\ifx\@\empty\let\next=\relax
   \else\let\next=\setfl@gs\expandafter\xdef
   \csname#1tr@cetrue\endcsname{}\fi\next}
\def\m@ketag#1#2{\expandafter\n@wcount\csname#2tagno\endcsname
     \csname#2tagno\endcsname=0\let\tail=\all\xdef\all{\tail#2,}
   \ifx#1\l@c@l\let\tail=\r@s@t\xdef\r@s@t{\csname#2tagno\endcsname=0\tail}\fi
   \expandafter\gdef\csname#2cite\endcsname##1{\expandafter
     \ifx\csname#2tag##1\endcsname\relax?\else\csname#2tag##1\endcsname\fi
     \expandafter\ifx\csname#2tr@cetrue\endcsname\relax\else
     \write\cit@tionsout{#2tag ##1 cited on page \folio.}\fi}
   \expandafter\gdef\csname#2page\endcsname##1{\expandafter
     \ifx\csname#2page##1\endcsname\relax?\else\csname#2page##1\endcsname\fi
     \expandafter\ifx\csname#2tr@cetrue\endcsname\relax\else
     \write\cit@tionsout{#2tag ##1 cited on page \folio.}\fi}
   \expandafter\gdef\csname#2tag\endcsname##1{\expandafter
      \ifx\csname#2check##1\endcsname\relax
      \expandafter\xdef\csname#2check##1\endcsname{}%
      \else\immediate\write16{Warning: #2tag ##1 used more than once.}\fi
      \multit@g{#1}{#2}##1/X%
      \write\t@gsout{#2tag ##1 assigned number \csname#2tag##1\endcsname\space
      on page \number\count0.}%
   \csname#2tag##1\endcsname}}

\def\multit@g#1#2#3/#4X{\def\t@mp{#4}\ifx\t@mp\empty%
      \global\advance\csname#2tagno\endcsname by 1 
      \expandafter\xdef\csname#2tag#3\endcsname
      {#1\number\csname#2tagno\endcsnameX}%
   \else\expandafter\ifx\csname#2last#3\endcsname\relax
      \expandafter\n@wcount\csname#2last#3\endcsname
      \global\advance\csname#2tagno\endcsname by 1 
      \expandafter\xdef\csname#2tag#3\endcsname
      {#1\number\csname#2tagno\endcsnameX}
      \write\t@gsout{#2tag #3 assigned number \csname#2tag#3\endcsname\space
      on page \number\count0.}\fi
   \global\advance\csname#2last#3\endcsname by 1
   \def\t@mp{\expandafter\xdef\csname#2tag#3/}%
   \expandafter\t@mp\@mputate#4\endcsname
   {\csname#2tag#3\endcsname\lastpart{\csname#2last#3\endcsname}}\fi}
\def\t@gs#1{\def\all{}\m@ketag#1e\m@ketag#1s\m@ketag\t@d@l p
\let\realscite\scite
\let\realstag\stag
   \m@ketag\gl@b@l r \n@wread\t@gsin
   \openin\t@gsin=\jobname.tgs \re@der \closein\t@gsin
   \n@wwrite\t@gsout\openout\t@gsout=\jobname.tgs }
\outer\def\localtags{\t@gs\l@c@l}
\outer\def\globaltags{\t@gs\gl@b@l}
\outer\def\newlocaltag#1{\m@ketag\l@c@l{#1}}
\outer\def\newglobaltag#1{\m@ketag\gl@b@l{#1}}

\newif\ifpr@ 
\def\m@kecs #1tag #2 assigned number #3 on page #4.%
   {\expandafter\gdef\csname#1tag#2\endcsname{#3}
   \expandafter\gdef\csname#1page#2\endcsname{#4}
   \ifpr@\expandafter\xdef\csname#1check#2\endcsname{}\fi}
\def\re@der{\ifeof\t@gsin\let\next=\relax\else
   \read\t@gsin to\t@gline\ifx\t@gline\v@idline\else
   \expandafter\m@kecs \t@gline\fi\let \next=\re@der\fi\next}
\def\pretags#1{\pr@true\pret@gs#1,,}
\def\pret@gs#1,{\def\@{#1}\ifx\@\empty\let\n@xtfile=\relax
   \else\let\n@xtfile=\pret@gs \openin\t@gsin=#1.tgs \message{#1} \re@der 
   \closein\t@gsin\fi \n@xtfile}

\newcount\sectno\sectno=0\newcount\subsectno\subsectno=0
\newif\ifultr@local \def\ultralocal{\ultr@localtrue}
\def\firstpart{\number\sectno}
\def\lastpart#1{\ifcase#1 \or a\or b\or c\or d\or e\or f\or g\or h\or 
   i\or k\or l\or m\or n\or o\or p\or q\or r\or s\or t\or u\or v\or w\or 
   x\or y\or z \fi}

\def\resetall{\global\advance\sectno by 1\subsectno=0
   \gdef\firstpart{\number\sectno}\r@s@t}
\def\resetsub{\global\advance\subsectno by 1
   \gdef\firstpart{\number\sectno.\number\subsectno}\r@s@t}
\def\newsection#1\par{\resetall\vskip0pt plus.3\vsize\penalty-250
   \vskip0pt plus-.3\vsize\bigskip\bigskip
   \message{#1}\leftline{\bf#1}\nobreak\bigskip}
\def\subsection#1\par{\ifultr@local\resetsub\fi
   \vskip0pt plus.2\vsize\penalty-250\vskip0pt plus-.2\vsize
   \bigskip\smallskip\message{#1}\leftline{\bf#1}\nobreak\medskip}


\newdimen\marginshift

\newdimen\margindelta
\newdimen\marginmax
\newdimen\marginmin

\def\margininit{       
\marginmax=3 true cm                  
				      
\margindelta=0.1 true cm              
\marginmin=0.1true cm                 
\marginshift=\marginmin
}    

\def\t@gsjj#1,{\def\@{#1}\ifx\@\empty\let\next=\relax\else\let\next=\t@gsjj
   \def\@@{p}\ifx\@\@@\else
   \expandafter\gdef\csname#1cite\endcsname##1{\citejj{##1}}
   \expandafter\gdef\csname#1page\endcsname##1{?}
   \expandafter\gdef\csname#1tag\endcsname##1{\tagjj{##1}}\fi\fi\next}
\newif\ifshowstuffinmargin
\showstuffinmarginfalse
\def\jjtags{\showstuffinmargintrue
\ifx\all\relax\else\expandafter\t@gsjj\all,\fi}

\def\tagjj#1{\realstag{#1}\mginpar{\zeigen{#1}}}
\def\citejj#1{\zeigen{#1}\mginpar{\rechnen{#1}}}

\def\rechnen#1{\expandafter\ifx\csname stag#1\endcsname\relax ??\else
                           \csname stag#1\endcsname\fi}

\newdimen\theight

\def\marginfont{\sevenrm}

\def\trymarginbox#1{\setbox0=\hbox{\marginfont\hskip\marginshift #1}%
		\global\marginshift\wd0 
		\global\advance\marginshift\margindelta}

\def \mginpar#1{%
\ifvmode\setbox0\hbox to \hsize{\hfill\rlap{\marginfont\quad#1}}%
\ht0 0cm
\dp0 0cm
\box0\vskip-\baselineskip
\else 
             \vadjust{\trymarginbox{#1}%
		\ifdim\marginshift>\marginmax \global\marginshift\marginmin
			\trymarginbox{#1}%
                \fi
             \theight=\ht0
             \advance\theight by \dp0    \advance\theight by \lineskip
             \kern -\theight \vbox to \theight{\rightline{\rlap{\box0}}%
\vss}}\fi}


\def\t@gsoff#1,{\def\@{#1}\ifx\@\empty\let\next=\relax\else\let\next=\t@gsoff
   \def\@@{p}\ifx\@\@@\else
   \expandafter\gdef\csname#1cite\endcsname##1{\zeigen{##1}}
   \expandafter\gdef\csname#1page\endcsname##1{?}
   \expandafter\gdef\csname#1tag\endcsname##1{\zeigen{##1}}\fi\fi\next}
\def\verbatimtags{\showstuffinmarginfalse
\ifx\all\relax\else\expandafter\t@gsoff\all,\fi}
\def\zeigen#1{\hbox{$\langle$}#1\hbox{$\rangle$}}

\def\(#1){\edef\dot@g{\ifmmode\ifinner(\hbox{\noexpand\etag{#1}})
   \else\noexpand\eqno(\hbox{\noexpand\etag{#1}})\fi
   \else(\noexpand\ecite{#1})\fi}\dot@g}

\newif\ifbr@ck
\def\eat#1{}
\def\[#1]{\br@cktrue[\br@cket#1'X]}
\def\br@cket#1'#2X{\def\temp{#2}\ifx\temp\empty\let\next\eat
   \else\let\next\br@cket\fi
   \ifbr@ck\br@ckfalse\br@ck@t#1,X\else\br@cktrue#1\fi\next#2X}
\def\br@ck@t#1,#2X{\def\temp{#2}\ifx\temp\empty\let\neext\eat
   \else\let\neext\br@ck@t\def\temp{,}\fi
   \def\teemp{#1}\ifx\teemp\empty\else\rcite{#1}\fi\temp\neext#2X}
\def\resetbr@cket{\gdef\[##1]{[\rtag{##1}]}}
\def\references{\resetbr@cket\newsection References\par}

\newtoks\symb@ls\newtoks\s@mb@ls\newtoks\p@gelist\n@wcount\ftn@mber
    \ftn@mber=1\newif\ifftn@mbers\ftn@mbersfalse\newif\ifbyp@ge\byp@gefalse
\def\defm@rk{\ifftn@mbers\n@mberm@rk\else\symb@lm@rk\fi}
\def\n@mberm@rk{\xdef\m@rk{{\the\ftn@mber}}%
    \global\advance\ftn@mber by 1 }
\def\rot@te#1{\let\temp=#1\global#1=\expandafter\r@t@te\the\temp,X}
\def\r@t@te#1,#2X{{#2#1}\xdef\m@rk{{#1}}}
\def\b@@st#1{{$^{#1}$}}\def\str@p#1{#1}
\def\symb@lm@rk{\ifbyp@ge\rot@te\p@gelist\ifnum\expandafter\str@p\m@rk=1 
    \s@mb@ls=\symb@ls\fi\write\f@nsout{\number\count0}\fi \rot@te\s@mb@ls}
\def\byp@ge{\byp@getrue\n@wwrite\f@nsin\openin\f@nsin=\jobname.fns 
    \n@wcount\currentp@ge\currentp@ge=0\p@gelist={0}
    \re@dfns\closein\f@nsin\rot@te\p@gelist
    \n@wread\f@nsout\openout\f@nsout=\jobname.fns }
\def\m@kelist#1X#2{{#1,#2}}
\def\re@dfns{\ifeof\f@nsin\let\next=\relax\else\read\f@nsin to \f@nline
    \ifx\f@nline\v@idline\else\let\t@mplist=\p@gelist
    \ifnum\currentp@ge=\f@nline
    \global\p@gelist=\expandafter\m@kelist\the\t@mplistX0
    \else\currentp@ge=\f@nline
    \global\p@gelist=\expandafter\m@kelist\the\t@mplistX1\fi\fi
    \let\next=\re@dfns\fi\next}
\def\symbols#1{\symb@ls={#1}\s@mb@ls=\symb@ls} 
\def\bigsymbol{\textstyle}
\symbols{\bigsymbol\ast,\dagger,\ddagger,\sharp,\flat,\natural,\star}
\def\ftnumbers{\ftn@mberstrue} \def\ftsymbols{\ftn@mbersfalse}
\def\paginal{\byp@ge} \def\resetftnumbers{\ftn@mber=1}
\def\ftnote#1{\defm@rk\expandafter\expandafter\expandafter\footnote
    \expandafter\b@@st\m@rk{#1}}

\long\def\jump#1\endjump{}
\def\ssum{\mathop{\lower .1em\hbox{$\textstyle\Sigma$}}\nolimits}

\def\qed{\nobreak\kern 1em \vrule height .5em width .5em depth 0em}
\def\newneq{\hbox{\rlap{\hbox to 1\wd9{\hss$=$\hss}}\raise .1em 
   \hbox to 1\wd9{\hss$\scriptscriptstyle/$\hss}}}
\def\subsetne{\setbox9 = \hbox{$\subset$}\mathrel{\hbox{\rlap
   {\lower .4em \newneq}\raise .13em \hbox{$\subset$}}}}
\def\supsetne{\setbox9 = \hbox{$\subset$}\mathrel{\hbox{\rlap
   {\lower .4em \newneq}\raise .13em \hbox{$\supset$}}}}

\def\vbar{\mathchoice{\vrule height6.3ptdepth-.5ptwidth.8pt\kern-.8pt}
   {\vrule height6.3ptdepth-.5ptwidth.8pt\kern-.8pt}
   {\vrule height4.1ptdepth-.35ptwidth.6pt\kern-.6pt}
   {\vrule height3.1ptdepth-.25ptwidth.5pt\kern-.5pt}}
\def\f@dge{\mathchoice{}{}{\mkern.5mu}{\mkern.8mu}}
\def\b@c#1#2{{\rm \mkern#2mu\vbar\mkern-#2mu#1}}
\def\b@b#1{{\rm I\mkern-3.5mu #1}}
\def\b@a#1#2{{\rm #1\mkern-#2mu\f@dge #1}}
\def\bb#1{{\count4=`#1 \advance\count4by-64 \ifcase\count4\or\b@a A{11.5}\or
   \b@b B\or\b@c C{5}\or\b@b D\or\b@b E\or\b@b F \or\b@c G{5}\or\b@b H\or
   \b@b I\or\b@c J{3}\or\b@b K\or\b@b L \or\b@b M\or\b@b N\or\b@c O{5} \or
   \b@b P\or\b@c Q{5}\or\b@b R\or\b@a S{8}\or\b@a T{10.5}\or\b@c U{5}\or
   \b@a V{12}\or\b@a W{16.5}\or\b@a X{11}\or\b@a Y{11.7}\or\b@a Z{7.5}\fi}}

\catcode`\X=11 \catcode`\@=12


\expandafter\ifx\csname citeadd.tex\endcsname\relax
\expandafter\gdef\csname citeadd.tex\endcsname{}
\else \message{Hey!  Apparently you were trying to
\string\input{citeadd.tex} twice.   This does not make sense.} 
\errmessage{Please edit your file (probably \jobname.tex) and remove
any duplicate ``\string\input'' lines}\endinput\fi

\sectno=-1   
\localtags
\NoBlackBoxes
\define\mr{\medskip\roster}
\define\sn{\smallskip\noindent}
\define\mn{\medskip\noindent}
\define\bn{\bigskip\noindent}
\define\ub{\underbar}
\define\wilog{\text{without loss of generality}}
\define\ermn{\endroster\medskip\noindent}

\define\dbcu{\dsize\bigcup}
\define \nl{\newline}
\magnification=\magstep 1
\documentstyle {amsppt}
\pageheight{8.5truein}
\topmatter
\title{Classification theory for theories with NIP - a modest beginning \\
 Sh715} \endtitle
\rightheadtext{CTT with NIP}
\author {Saharon Shelah \thanks {\null\newline I would like to thank 
Alice Leonhardt for the beautiful typing. \null\newline
  The author would like to thank the United States Israel Binational Science
Foundation for partially supporting this research. \null\newline 
  Latest Revision - 00/Feb/28 \null\newline
  Done 10/98, 9/99} \endthanks} \endauthor 
\affil{Institute of Mathematics\\
 The Hebrew University\\
 Jerusalem, Israel
 \medskip
 Rutgers University\\
 Mathematics Department\\
 New Brunswick, NJ  USA} \endaffil

\abstract  A relevant thesis is that for the family of complete first
order theories with NIP (i.e. without the independence property) there
is a substantial theory, like the family of stable (and the family of
simple) first order theories.  We examine some properties.
\endabstract
\endtopmatter
\document  

\expandafter\ifx\csname alice2jlem.tex\endcsname\relax
  \expandafter\xdef\csname alice2jlem.tex\endcsname{\the\catcode`@}
\else \message{Hey!  Apparently you were trying to
\string\input{alice2jlem.tex}  twice.   This does not make sense.}
\errmessage{Please edit your file (probably \jobname.tex) and remove
any duplicate ``\string\input'' lines}\endinput\fi

\expandafter\ifx\csname bib4plain.tex\endcsname\relax
  \expandafter\gdef\csname bib4plain.tex\endcsname{}
\else \message{Hey!  Apparently you were trying to \string\input
  bib4plain.tex twice.   This does not make sense.}
\errmessage{Please edit your file (probably \jobname.tex) and remove
any duplicate ``\string\input'' lines}\endinput\fi

\def\renewcommand{\newcommand}	       
\edef\cite{\the\catcode`@}%
\catcode`@ = 11
\let\@oldatcatcode = \cite
\chardef\@letter = 11
\chardef\@other = 12
%
%
%
%
\def\@innerdef#1#2{\edef#1{\expandafter\noexpand\csname #2\endcsname}}%
%
%
\@innerdef\@innernewcount{newcount}%
\@innerdef\@innernewdimen{newdimen}%
\@innerdef\@innernewif{newif}%
\@innerdef\@innernewwrite{newwrite}%
%
%
%
\def\@gobble#1{}%
%
%
%
\ifx\inputlineno\@undefined
   \let\@linenumber = \empty 
\else
   \def\@linenumber{\the\inputlineno:\space}%
\fi
%
%
%
\def\@futurenonspacelet#1{\def\cs{#1}%
   \afterassignment\@stepone\let\@nexttoken=
}%
\begingroup 
\def\\{\global\let\@stoken= }%
\\ 
\endgroup
\def\@stepone{\expandafter\futurelet\cs\@steptwo}%
\def\@steptwo{\expandafter\ifx\cs\@stoken\let\@@next=\@stepthree
   \else\let\@@next=\@nexttoken\fi \@@next}%
\def\@stepthree{\afterassignment\@stepone\let\@@next= }%
%
%
%
\def\@getoptionalarg#1{%
   \let\@optionaltemp = #1%
   \let\@optionalnext = \relax
   \@futurenonspacelet\@optionalnext\@bracketcheck
}%
%
%
\def\@bracketcheck{%
   \ifx [\@optionalnext
      \expandafter\@@getoptionalarg
   \else
      \let\@optionalarg = \empty
      \expandafter\@optionaltemp
   \fi
}%
\def\@@getoptionalarg[#1]{%
   \def\@optionalarg{#1}%
   \@optionaltemp
}%
%
%
%
\def\@nnil{\@nil}%
\def\@fornoop#1\@@#2#3{}%
\def\@for#1:=#2\do#3{%
   \edef\@fortmp{#2}%
   \ifx\@fortmp\empty \else
      \expandafter\@forloop#2,\@nil,\@nil\@@#1{#3}%
   \fi
}%
\def\@forloop#1,#2,#3\@@#4#5{\def#4{#1}\ifx #4\@nnil \else
       #5\def#4{#2}\ifx #4\@nnil \else#5\@iforloop #3\@@#4{#5}\fi\fi
}%
\def\@iforloop#1,#2\@@#3#4{\def#3{#1}\ifx #3\@nnil
       \let\@nextwhile=\@fornoop \else
      #4\relax\let\@nextwhile=\@iforloop\fi\@nextwhile#2\@@#3{#4}%
}%
%
%
%
\@innernewif\if@fileexists
\def\@testfileexistence{\@getoptionalarg\@finishtestfileexistence}%
\def\@finishtestfileexistence#1{%
   \begingroup
      \def\extension{#1}%
      \immediate\openin0 =
         \ifx\@optionalarg\empty\jobname\else\@optionalarg\fi
         \ifx\extension\empty \else .#1\fi
         \space
      \ifeof 0
         \global\@fileexistsfalse
      \else
         \global\@fileexiststrue
      \fi
      \immediate\closein0
   \endgroup
}%
%
%
%
%
\def\bibliographystyle#1{%
   \@readauxfile
   \@writeaux{\string\bibstyle{#1}}%
}%
\let\bibstyle = \@gobble
%
%
\let\bblfilebasename = \jobname
\def\bibliography#1{%
   \@readauxfile
   \@writeaux{\string\bibdata{#1}}%
   \@testfileexistence[\bblfilebasename]{bbl}%
   \if@fileexists
      \nobreak
      \@readbblfile
   \fi
}%
\let\bibdata = \@gobble
%
%
\def\nocite#1{%
   \@readauxfile
   \@writeaux{\string\citation{#1}}%
}%
\@innernewif\if@notfirstcitation
%
%
\def\cite{\@getoptionalarg\@cite}%
%
%
\def\@cite#1{%
   \let\@citenotetext = \@optionalarg
   \printcitestart
   \nocite{#1}%
   \@notfirstcitationfalse
   \@for \@citation :=#1\do
   {%
      \expandafter\@onecitation\@citation\@@
   }%
   \ifx\empty\@citenotetext\else
      \printcitenote{\@citenotetext}%
   \fi
   \printcitefinish
}%
\def\@onecitation#1\@@{%
   \if@notfirstcitation
      \printbetweencitations
   \fi
   \expandafter \ifx \csname\@citelabel{#1}\endcsname \relax
      \if@citewarning
         \message{\@linenumber Undefined citation `#1'.}%
      \fi
      \expandafter\gdef\csname\@citelabel{#1}\endcsname{%
\strut
\vadjust{\vskip-\dp\strutbox
\vbox to 0pt{\vss\parindent0cm \leftskip=\hsize 
\advance\leftskip3mm
\advance\hsize 4cm\strut\openup-4pt 
\rightskip 0cm plus 1cm minus 0.5cm ?  #1 ?\strut}}
         {\tt
            \escapechar = -1
            \nobreak\hskip0pt
            \expandafter\string\csname#1\endcsname
            \nobreak\hskip0pt
         }%
      }%
   \fi
   \csname\@citelabel{#1}\endcsname
   \@notfirstcitationtrue
}%
%
%
\def\@citelabel#1{b@#1}%
%
%
\def\@citedef#1#2{\expandafter\gdef\csname\@citelabel{#1}\endcsname{#2}}%
%
%
%
\def\@readbblfile{%
   \ifx\@itemnum\@undefined
      \@innernewcount\@itemnum
   \fi
   \begingroup
      \def\begin##1##2{%
         \setbox0 = \hbox{\biblabelcontents{##2}}%
         \biblabelwidth = \wd0
      }%
      \def\end##1{}
      %
      %
      \@itemnum = 0
      \def\bibitem{\@getoptionalarg\@bibitem}%
      \def\@bibitem{%
         \ifx\@optionalarg\empty
            \expandafter\@numberedbibitem
         \else
            \expandafter\@alphabibitem
         \fi
      }%
      \def\@alphabibitem##1{%
         \expandafter \xdef\csname\@citelabel{##1}\endcsname {\@optionalarg}%
         \ifx\biblabelprecontents\@undefined
            \let\biblabelprecontents = \relax
         \fi
         \ifx\biblabelpostcontents\@undefined
            \let\biblabelpostcontents = \hss
         \fi
         \@finishbibitem{##1}%
      }%
      \def\@numberedbibitem##1{%
         \advance\@itemnum by 1
         \expandafter \xdef\csname\@citelabel{##1}\endcsname{\number\@itemnum}%
         \ifx\biblabelprecontents\@undefined
            \let\biblabelprecontents = \hss
         \fi
         \ifx\biblabelpostcontents\@undefined
            \let\biblabelpostcontents = \relax
         \fi
         \@finishbibitem{##1}%
      }%
      \def\@finishbibitem##1{%
         \biblabelprint{\csname\@citelabel{##1}\endcsname}%
         \@writeaux{\string\@citedef{##1}{\csname\@citelabel{##1}\endcsname}}%
         \ignorespaces
      }%
      %
      %
      \let\em = \bblem
      \let\newblock = \bblnewblock
      \let\sc = \bblsc
      \frenchspacing
      \clubpenalty = 4000 \widowpenalty = 4000
      \tolerance = 10000 \hfuzz = .5pt
      \everypar = {\hangindent = \biblabelwidth
                      \advance\hangindent by \biblabelextraspace}%
      \bblrm
      \parskip = 1.5ex plus .5ex minus .5ex
      \biblabelextraspace = .5em
      \bblhook
      \input \bblfilebasename.bbl
   \endgroup
}%
%
%
\@innernewdimen\biblabelwidth
\@innernewdimen\biblabelextraspace
%
%
%
\def\biblabelprint#1{%
   \noindent
   \hbox to \biblabelwidth{%
      \biblabelprecontents
      \biblabelcontents{#1}%
      \biblabelpostcontents
   }%
   \kern\biblabelextraspace
}%
%
%
%
\def\biblabelcontents#1{{\bblrm [#1]}}%
%
%
\def\bblrm{\rm}%
%
%
\def\bblem{\it}%
%
%
\def\bblsc{\ifx\@scfont\@undefined
              \font\@scfont = cmcsc10
           \fi
           \@scfont
}%
%
%
\def\bblnewblock{\hskip .11em plus .33em minus .07em }%
%
%
\let\bblhook = \empty
%
%
%
\def\printcitestart{[}
\def\printcitefinish{]}
\def\printbetweencitations{, }
\def\printcitenote#1{, #1}
%
%
%
\let\citation = \@gobble
%
%
%
\@innernewcount\@numparams
%
%
\def\newcommand#1{%
   \def\@commandname{#1}%
   \@getoptionalarg\@continuenewcommand
}%
%
%
\def\@continuenewcommand{%
   \@numparams = \ifx\@optionalarg\empty 0\else\@optionalarg \fi \relax
   \@newcommand
}%
%
%
\def\@newcommand#1{%
   \def\@startdef{\expandafter\edef\@commandname}%
   \ifnum\@numparams=0
      \let\@paramdef = \empty
   \else
      \ifnum\@numparams>9
         \errmessage{\the\@numparams\space is too many parameters}%
      \else
         \ifnum\@numparams<0
            \errmessage{\the\@numparams\space is too few parameters}%
         \else
            \edef\@paramdef{%
               \ifcase\@numparams
                  \empty  No arguments.
               \or ####1%
               \or ####1####2%
               \or ####1####2####3%
               \or ####1####2####3####4%
               \or ####1####2####3####4####5%
               \or ####1####2####3####4####5####6%
               \or ####1####2####3####4####5####6####7%
               \or ####1####2####3####4####5####6####7####8%
               \or ####1####2####3####4####5####6####7####8####9%
               \fi
            }%
         \fi
      \fi
   \fi
   \expandafter\@startdef\@paramdef{#1}%
}%
%
%
%
%
\def\@readauxfile{%
   \if@auxfiledone \else 
      \global\@auxfiledonetrue
      \@testfileexistence{aux}%
      \if@fileexists
         \begingroup
            \endlinechar = -1
            \catcode`@ = 11
            \input \jobname.aux
         \endgroup
      \else
         \message{\@undefinedmessage}%
         \global\@citewarningfalse
      \fi
      \immediate\openout\@auxfile = \jobname.aux
   \fi
}%
%
%
\newif\if@auxfiledone
\ifx\noauxfile\@undefined \else \@auxfiledonetrue\fi
%
%
%
%
\@innernewwrite\@auxfile
\def\@writeaux#1{\ifx\noauxfile\@undefined \write\@auxfile{#1}\fi}%
%
%
%
\ifx\@undefinedmessage\@undefined
   \def\@undefinedmessage{No .aux file; I won't give you warnings about
                          undefined citations.}%
\fi
%
%
\@innernewif\if@citewarning
\ifx\noauxfile\@undefined \@citewarningtrue\fi
%
%
%
\catcode`@ = \@oldatcatcode


\def\widestnumber#1#2{}

\def\rm{\fam0 \tenrm}

\def\fakesubhead#1\endsubhead{\bigskip\noindent{\bf#1}\par}



%
%
%

%

\font\textrsfs=rsfs10
\font\scriptrsfs=rsfs7
\font\scriptscriptrsfs=rsfs5

\newfam\rsfsfam
\textfont\rsfsfam=\textrsfs
\scriptfont\rsfsfam=\scriptrsfs
\scriptscriptfont\rsfsfam=\scriptscriptrsfs

\edef\oldcatcodeofat{\the\catcode`\@}
\catcode`\@11

\def\Cal@@#1{\noaccents@ \fam \rsfsfam #1}

\catcode`\@\oldcatcodeofat


\expandafter\ifx \csname margininit\endcsname \relax\else\margininit\fi

\newpage

\head {Anotated Content} \endhead  \resetall 
\bn
\S1 Indiscernible sequences and averages
\mr
\item "{{}}"  [We consider indiscernible sequence $\bold{\bar b} = \langle
\bar b_t:t \in I \rangle$ wondering, do they have an average as in the stable
case.  We investigate the set of $\varphi(\bar x,\bar y)$ such that every
instance $\varphi(\bar x,\bar a)$ divide $\bold{\bar b}$ to a finite/co-finite
sets, and some can divide it only to finitely many intervals; this is always
the case if $T$ has NIP.  If $T$ has NIP, indiscernible sequences
behave reasonably while indiscernible sets behave
nicely and so does $p \in S(M)$ connected with them which we call
stable types.  We then investigate
having an unstable/nip $\varphi(x,y;\bar c)$, i.e. on singletons.]
\endroster
\bn
\S2 Characteristics of types
\mr
\item "{{}}"  [Each indiscernible sequence $\bold{\bar b} = \langle \bar b_t:
t \in I \rangle$, has for each $\varphi = \varphi(\bar x,\bar y)$ a characteristic
number $n = n_{\bold{\bar b},\varphi}$, the maximal number of intervals to
which an instance $\varphi(\bar x,\bar b)$ can divide $\bold{\bar b}$.  We
wonder what we can say about it.]
\endroster
\bn
\S3  Shrinking indiscernbles
\mr
\item "{{}}" [For an indiscernible sequence to a set, if we increase the
set a little, not much indiscernability is lost.]
\endroster
\bn
\S4  Perpendicular endless indiscernible sequences
\mr
\item "{{}}"  [We define perpendicularity and investigate its basic property;
any two mutually indiscernible sequences are perpendicular.  E.g. (for NIP
theories) one sequence cannot be non perpendicular to $\ge |T|^+$ pairwise
perpendicular sequences.  We then deal with 
$\bold F^{\text{sp}}_{|T|^+}$-constructions.]
\endroster
\bn
\S5 Indiscernible sequences perpendicular to cuts
\mr
\item "{{}}"  [Using construction as above we show that we can build models
controlling quite tightly the dual cofinality of such sequences.]
\endroster
\newpage

\head {\S1 Indiscernible sequences and averages} \endhead  \resetall \sectno=1
\bn
We try to continue \cite[Ch.II,4.13]{Sh:c}, but we do not rely on it.
NIP stands for = no independence property.  For stable (complete first
order) theories, the notions of indiscernible set and its average (and
local versions of them) play important role.  For unstable theory,
indiscernible sequences are not necessarily indiscernible sets.  Still
for indiscernible sets $\bold I$ if $T$ has NIP, the basic claim
guaranteeing the existence of averages (any $\varphi(\bar x,\bar b)$
divide $\bold I$ to a finite and co-finite set) hold, and for
indiscernible sequences the division is into the union of
$<_{n_\varphi}$ convex sets.  For any $T$, we can still look at the
first order formulas $\varphi(\bar x,\bar y)$ which behaves well,
i.e. any $\varphi(\bar x,\bar b)$ divide any appropriate $\bold I$ as
above.

In \scite{np1.1} - \scite{np1.4} + \scite{np1.4}(c) + (d) we define
the relevant notions: average (Av$(\bold J,D)$ or Av$_\varphi(\bold
J,D)$ or Av$(\bold J,\langle \bar b_t:t \in I \rangle)$ and av for
finite sequence averaging formulas for $\bold I$(avf,daf), and state
some basic properties.

In particular we look at indiscernible sequence of ``finite distance",
those are related to canonical bases (of types, of indiscernible sets)
play important role for stable theories, hence we try to define
parallels in \scite{np1.4}, see \scite{np1.4c}(2).

Next we note a dichotomy for the types $p \in S^m(M)$.  Such a type
$p$ may be stable (see Definition \scite{np1.8}, Claim \scite{np1.5} -
\scite{np1.6}); not only is the type definable, but for every
ultrafilter $D$ on ${}^m M$ with Av$(M,D) = p$, any indiscernible set
constructed from $D$ is an indiscernible set, and the definition comes
from appropriate finite large enough $(\Delta,k)$-indiscernible sets.
If $p \in S^m(M)$ is not stable, then there is a partial order with
infinite chains closely related to it.  We conclude that if $T$ is
unstable (with NIP), then some $\varphi(x,y,\bar c)$ define a quasi
order with infinite chains and if $T$ is unstable some
$\varphi(x,y;\bar c)$ has the order property (though not necessarily
the property (E) of Eherenfeucht).
\bigskip

\demo{\stag{np1.0} Context}  $T$ a complete first order theory, its monster
model being ${\frak C} = {\frak C}_T$ as usual in \cite{Sh:c}.
\enddemo
\bigskip

\definition{\stag{np1.0a} Definition}  
$T$ is NIP means it does not have the independence property, IP in short,
i.e. for no $\varphi(\bar x,\bar y)$ do we have for every $n$
\mr
\item "{$\boxtimes^n_\varphi$}"  ${\frak C} \models (\exists \bar y_0,\dotsc,
\bar y_{n-1}) \dsize \bigwedge_{\eta \in {}^n 2} (\exists \bar x)
(\dsize \bigwedge_{\ell < n} 
\varphi(\bar x,\bar y_\ell)^{\text{if}(\eta(\ell))}$).
\endroster
\enddefinition
\bigskip

\definition{\stag{np1.1} Definition}  1)  For $\bold J \subseteq
{}^{\omega >} {\frak C},m < \omega$, a set $\bold I \subseteq {}^m{\frak C}$, 
an ultrafilter on $D$ over $\bold I$ we let Av$(\bold J,D)$
be $\{\varphi(\bar x,\bar a):\bar x = \langle x_\ell:\ell < m \rangle,
\bar a \in \bold J$ and $\{\bar b \in \bold I:\models \varphi
(\bar b,\bar a)\} \in D\}$.  It will be called the $D$-average over 
$\bold J$.  If $\bold J = {}^{\omega >}B$ we may write $B$ instead of 
$\bold J$ (or $M$ if $B = |M|$).  (Av stands for average). \nl
2) If $D$ is an ultrafilter over $\bold I \subseteq {}^m{\frak C},I$ an 
infinite linear order, $\bar{\bold b} = \langle \bar b_t:t \in I \rangle$, 
we say that $\bar{\bold b}$ is an 
$(A,D)$-indiscernible sequence \ub{if} for each $t \in I$ we have
tp$(\bar b_t,A \cup \{\bar b_s:s <_I t\}) =$ Av$(A \cup \{\bar b_s:
s <_I t\},D)$.  If $A = \cup \{\bar c:\bar c \in \text{ Dom}(D)\}$, we may 
write ``$\bar{\bold b}$ is a $D$-indiscernible sequence". \nl
3) Av$_\varphi(A,D)$ where $\varphi = \varphi(\bar x,\bar y)$ is
$\{\varphi(\bar x,\bar a):\varphi(\bar x,\bar a) \in$ Av$(A,D)\}$ and
Av$_\Delta(A,D) = \dbcu_{\varphi \in \Delta}$ Av$_\varphi(A,D)$. 
\enddefinition
\bigskip

\proclaim{\stag{np1.1a} Claim}  1) For $D$ an ultrafilter on $\bold I
\subseteq {}^m {\frak C}$ and $B \subseteq {\frak C}$ we have {\rm Av\/}$(B,D) \in
S^m(B)$. \nl
2) If $\bold I \subseteq {}^m A$ and $D$ is an ultrafilter on $\bold I$ and
$\langle \bar b_t:t \in I \rangle$ is an $(A,D)$-indiscernible sequence
\ub{then} $\langle \bar b_t:t \in I \rangle$ is an indiscernible sequence
over $A$.
\endproclaim
\bigskip

\demo{Proof}  Check.
\enddemo
\bigskip

\definition{\stag{np1.2} Definition}  1) For 
an infinite linear order $I$ and an
indiscernible sequence $\bar{\bold b} = \langle \bar b_t:t \in I \rangle$,
having $\ell g(\bar b_t) =m$ for $t \in I$, we define:
\mr
\item "{$(a)$}"  avf$_{\text{pa}}(\bar{\bold b}) = \{\varphi(\bar x,
\bar y,\bar c):\ell g(\bar y) = m$, and for every $\bar a \in 
{}^{\ell g(\bar x)}{\frak C}$, the set $\{t \in I:{\frak C} \models 
\varphi(\bar a,\bar b_t,\bar c)\}$ is finite
\ub{or} the set $\{t \in I:{\frak C} \models 
\neg \varphi(\bar b_t;\bar a,\bar c)\}$ is finite$\}$ \nl
(avf stands for averagable formulas, pa stands for parameters)
\sn
\item "{$(b)$}"  avf$(\bar{\bold b}) = \{\varphi(\bar x;\bar y):
\varphi(\bar x,\bar y) \in \text{ avf}_{\text{pa}}(\bar{\bold b})$, i.e. no 
parameters$\}$
\sn
\item "{$(c)$}"  daf$(\bar{\bold b}) = \{\varphi(\bar x,\bar y):\ell g
(\bar y) = m$ and for every $\bar a \in {}^{\ell g(\bar x)}{\frak C}$ the
set $\{t \in I:{\frak C} \models$ \nl

\hskip40pt  $\varphi[\bar a,\bar b_t]\}$ is a finite union of convex 
subsets of $I\}$. \nl
Let daf$_{pa}(\bar{\bold b})$ be defined similarly allowing parameters,
\sn
\item "{$(d)$}"  daf$^n(\bold{\bar b})$ when the union is of $\le n$
convex sets; similarly is the other cases.
\ermn
2) For a sequence $\bar{\bold b} = \langle \bar b_t:t < k \rangle$
with $\bar b_\ell \in {}^m {\frak C}$, and formula $\varphi = 
\varphi(\bar y,\bar z),\ell g (\bar y) = m$, we define

$$
\align
\text{av}_\varphi(A,\langle \bar b_\ell:\ell < k) \rangle) = 
\{\varphi(\bar y,\bar c)^{\bold t}:&\bold t \in \{\text{true, false}\}, \\
  &\bar c \in {}^{\ell g(\bar z)} A, \text{ and } |\{\ell:\models \varphi
(\bar b_\ell,\bar c)^{\bold t}\}| > k/2\}
\endalign
$$
\mn
(this is not necessarily a type, just a set of formulas). \nl
3) $E = E^k_{\varphi(\bar y,\bar z)}$, a formula in $L(\tau_T)$, written
$\bar z_1 E \bar z_2$ with $\ell g(\bar z_1) = \ell g(\bar z_2) = 
(\ell g(\bar x)) \times k$ (written $(\bar x_1,\ldots,\bar x_n)$ instead of
$\bar x_1 \char 94 \ldots \char 94 x_n$, abusing notation) is defined 
as follows: 
$(\bar x_0,\dotsc,\bar x_{k-1})E(\bar x'_0,\dotsc,x'_{k-1}) =:$ \nl
$(\forall \bar z)(\dsize \bigvee_{u \subseteq k,|u| > k/2} \, 
\dsize \bigwedge_{\ell \in u} \varphi(\bar x_\ell,\bar z)
\equiv \dsize \bigvee_{u \subseteq k,|u| > k/2} \, 
\dsize \bigwedge_{\ell \in u} \varphi(\bar x'_\ell,\bar z))$.  \nl
Of course, it is an equivalence relation.
\enddefinition
\bigskip

\proclaim{\stag{np1.3} Claim}  If $\bold{\bar b} = \langle \bar b_t:t \in I
\rangle$ is an infinite indiscernible sequence, 
$\ell g(\bar b_t) = m$ and $\varphi
(\bar y;\bar z) \in \text{ avf}(\bold{\bar b})$ so $\ell g(\bar z) = m$, 
\ub{then} for every $k$ large enough we have:
\mr
\item "{$(a)$}"  for any $\bar c$ of length $\ell g(\bar z)$, for some 
truth value $\bold t$ the set $\{t \in I:\models \varphi(\bar b_t,\bar
c)^{\bold t}\}$ 
has $< k/2$ members
\sn
\item "{$(b)$}"  if $t_0,\dotsc,t_{k-1}$ are distinct members of $I$
\ub{then} {\rm av\/}$_\varphi(\langle \bar b_{t_\ell}:\ell < k \rangle,{\frak C}) \in
S^{\ell g(\bar m)}_\varphi({\frak C})$, in fact for every nonprincipal
ultrafilter $D$ over $\{\bar b_t:t \in I\}$ and set $A$ we have {\rm av\/}$_\varphi
(\langle b_{t_\ell}:\ell < k \rangle,A)$ is a subset of {\rm Av\/}$(A,D)$, in fact is
Av$_\varphi(A,D)$
\sn
\item "{$(c)$}"  if $t_0,\dotsc,t_{k-1} \in I$ with no repetitions and
$s_0,\dotsc,s_{k-1} \in I$ with no repetition \ub{then} $(\bar b_{t_0},\dotsc,
\bar b_{t_{k-1}}) \, E^k_{\varphi(\bar x,\bar y)}(\bar b_{s_0},\dotsc,
b_{s_{k_1}})$
\sn
\item "{$(d)$}"  for some finite $\Delta$: if $I',I \subseteq J$ where $J$ 
is a linear order, $\bold{\bar b}' = \langle \bar b'_t:
t \in J \rangle,\bold{\bar b}' \restriction I = \bold{\bar b}$ and 
$\bold{\bar b}'$ is $\Delta$-indiscernible sequence,
$|I'| \ge k$, \ub{then}
{\roster
\itemitem{ $(\alpha)$ }   (b),(c) holds for $\bold{\bar b}' \restriction I'$
and 
\sn
\itemitem{ $(\beta)$ }  $\varphi(\bar y,\bar z) \in 
\text{ {\rm avf\/}}(\bar{\bold b}' \restriction I)$.
\endroster}
\endroster
\endproclaim
\bigskip

\demo{Proof} \nl

(a)  By compactness. \nl

(b),(c)  Just think of the definitions. \nl

(d)  By compactness.  \hfill$\square_{\scite{np1.3}}$
\enddemo
\bigskip

\definition{\stag{np1.4} Definition}  Let $\bold{\bar b} = \langle \bar b_t:
t \in I \rangle$ an infinite indiscernible sequence. \nl
1) We define (Cb stands for canonical bases):
\mr
\item "{$(a)$}"  for $\varphi(\bar y;\bar z) \in \text{ avf}
(\bold{\bar b})$ let 
$Cb_{\varphi(\bar y;\bar z)}(\bold{\bar b})$ be $(\bar b_{t_0},\dotsc,
\bar b_{t_{k-1}})/E^k_{\varphi(\bar y,\bar z)} \in {\frak C}^{\text{eq}}$,
with $k = k_{\varphi(\bar y,\bar z)}(\bold{\bar b})$ minimal as in
\scite{np1.3}(a)
\sn
\item "{$(b)$}"  $Cb(\bold{\bar b}) = \text{ dcl}\{Cb_{\varphi(\bar y,
\bar z)}(\bold{\bar b}):\varphi(\bar y,\bar z) \in \text{ avf}
(\bold{\bar b})\} \subseteq {\frak C}^{\text{eq}}$.
\ermn
2) If $I$ has no last element \ub{then} Av$(A,\bar{\bold b}) =
\{\varphi(\bar x,\bar a):\bar a \in {}^{\omega >}A$ and $\models \varphi
(\bar b_t,\bar a)$ for every large enough $t \in I\}$,
Av$_\varphi(A,\bar{\bold b}) = \{\varphi(\bar x,\bar a):$ 
for every large enough $t \in I$ we have ${\frak C}
\models \varphi(\bar b_t,\bar a)\}$.
\nl
3) Let Av$_{\text{avf}}(\bold{A,\bar b})$ be
Av$_\Delta(\bold{\bar b},A)$ for $\Delta = \text{ avf}(\bold{\bar b})$,
similarly for replacements to avf.
\enddefinition
\bigskip

\proclaim{\stag{np1.4a} Claim}  [$T$ has NIP]  Assume $I$ is a linear
order with no last element and $\bold{\bar b} = \langle \bar b_t:t \in I
\rangle$ an indiscernible sequence, $\ell g(\bar b_t) = n$. \nl
1) {\rm av\/}$_\varphi(A,\bold{\bar b}) \in S^m_\varphi(A)$, see Definition
\scite{np1.2}(a). \nl
2) {\rm Av\/}$(A,\bold{\bar b}) \in S^m(A)$, see Definition \scite{np1.4}(c).
\endproclaim
\bigskip

\demo{Proof}  By \cite[II.4.13]{Sh:c}.
\enddemo
\bigskip

To formalize clause (d) of \scite{np1.3} let
\definition{\stag{np1.4b} Definition}   1) For a set $\Delta$ of
formulas and $k \le \omega$ we say
that $\langle \bar b^1_t:t \in I_1 \rangle,\langle \bar b^2_t:t \in I_2 \rangle$
are immediate $(\Delta,k)$-nb-s if:
\mr
\item "{$(a)$}"  both are $\Delta$-indiscernible sequences of length
$\ge k$
\sn
\item "{$(b)$}"  for some $\Delta$-indiscernible sequence $\langle
\bar b_t:t \in I \rangle$ we have $I_\ell \subseteq I,(\forall t \in I_\ell)
\bar b^\ell_t= \bar b_t$ for $\ell =1,2$.
\ermn
2) The relation ``being $(\Delta,k)$-nb-s" is the closure of being an
``immediate $(\Delta,k)$-nb" to an equivalence relation.
We say ``of distance $k$" if there is a chain of immediate
$(\Delta,k$)-nb-s of length $k$ of length $\le k$ starting with one
ending in the other.  We write $\Delta$ instead of $(\Delta,\omega)$
if $\Delta = L(T)$ we may omit $\Delta$. \nl
3) If $\bar{\bold b}^1,\bar{\bold b}^2$ are infinite indiscernible
sequences, we say there are ``essentially nb-s" if for every finite
$\delta \in L(T),k < \omega$ they are $(\Delta,k)$-nb-s. \nl
4) If $\bar{\bold b}$ is an infinite indiscernible sequence over $A$
we let $C_A(\bar{\bold b}) = \{\bar b_t$: for some $\bar{\bold b}'$
essentially nb-s of $\bar{\bold b},\bar b$ appears in $\bar{\bold b}'\}$.
\enddefinition
\bigskip

\proclaim{\stag{np1.4c} Claim}  1) If $\langle \bar b_t:t \in I \rangle$ is an
infinite indiscernible sequence and $\varphi(\bar y,\bar z) \in
\text{ daf}^{\,n}(\langle \bar b_t:t \in I \rangle)$, \ub{then} for some finite
$\Delta$ and $k$, for any $(\Delta,k)$-nb-s $\langle \bar b'_t:t \in I' \rangle$
of $\langle \bar b_t:t \in I \rangle$ we have $\varphi(\bar y,\bar z) \in
\text{ daf}^{\,n}(\langle \bar b_t:t \in I' \rangle)$. \nl
2) The result in (1) holds also for {\rm avf\/}$^n(\langle \bar b_t:t
\in I \rangle)$.  
If $\bold{\bar b} = \langle \bar b_t:t \in I \rangle$ is an infinite
indiscernible sequence and $\varphi(\bar y,\bar z) \in$ {\rm avf\/}$(\bold{\bar b})$,
\ub{then} for some finite $\Delta$ and $k$ for any $(\Delta,k)$-nb
$\bold{\bar b}'$ of $\bold{\bar b} = \langle \bar b^0_t:t \in I \rangle$ we
have Cb$_{\varphi(\bar y,\bar z)}(\bold{\bar b}') = Cb_{\varphi(y,\bar z)}
(\bold{\bar b})$. \nl
3) If $T$ has NIP and $\bold{\bar b}$ is an infinite indiscernible sequence
\ub{then} every $\varphi(\bar y,\bar z)$ with the right length of $\bar y$ belongs
to daf$(\bold{\bar b})$ (see \cite[Ch.II,\S4]{Sh:c}). \nl
4) If $\bar b_1,\bar b_2$ are $(\Delta,k)$-nb-s of distance $n$ for
every finite $\Delta$ and $k < \omega$, \ub{then} they are
$(L(T),k)$-nb-s of distance $n$.
\endproclaim
\bigskip

\demo{Proof}  Easy.
\enddemo
\bn
\centerline{$* \qquad * \qquad *$}
\bigskip

\proclaim{\stag{np1.5} Definition/Claim}  Assume that
$M \prec {\frak C},p \in S^m(M)$ and for
$\ell = 1,2,D_\ell$ is an ultrafilter on ${}^m M$ and $\bold{\bar b}^\ell =
\langle \bar b^\ell_t:t \in I_\ell \rangle$ is an infinite 
$D_\ell$-indiscernible sequence over $M$ such that 
{\rm Av\/}$(M,D_\ell) =p$ (so does not
depend on $\ell$).  \ub{Then}
\mr
\item "{$(a)$}" for every finite set $\Delta_1$ of formulas, finite
$A \subseteq M$ and $j \in \{1,2\}$ there is a function ${\Cal F}$ 
into $D_j$ such that
{\roster
\itemitem{ $(*)$ }  if $\alpha \le \omega$ and for each $\ell < \alpha$
we have $\bar b_\ell \in {}^m M,
\bar b_\ell \in {\Cal F}(\bar b_0,\dotsc,\bar b_{\ell -1})(\in D_j)$
\ub{then} the sequence $\langle \bar b^j_t:t \in I \rangle \char 94
\langle \bar b_\ell:\ell \in \alpha^* \rangle$ (where the $*$ in $\alpha^*$ means 
invert the order) is $\Delta_1$-indiscernible sequence over $A$
\endroster}
\item "{$(b)$}"  for any $\varphi(\bar y,\bar z)$ for finite large enough
$\Delta \subseteq L(T)$, if $j \in \{1,2\}$ and $\bar b_\ell \in {}^m M$ for 
$\ell < k$ where $k < \omega$ is large enough are as in clause (a), 
we have: if $j \in \{1,2\}$ and $\varphi = \varphi(\bar y,
\bar z) \in \text{ {\rm avf\/}}(\bold{\bar b}^j),k$ also large enough
as in \scite{np1.3} \ub{then}
{\roster
\itemitem{ $(i)$ }  {\rm av\/}$_\varphi({\frak C},\bar{\bold b}^j) = 
\text{ {\rm av\/}}_\varphi({\frak C},\langle \bar b_\ell:\ell < k \rangle)$
\sn
\itemitem{ $(ii)$ }  $p \restriction \varphi(\bar y,\bar z)$ is definable
as $\bigl\{ \varphi(y,\bar c)^{\bold t}:\bar c \subseteq M 
\text{ and } |\{\ell < k:\neg \varphi (\bar b_\ell,\bar c)^{\bold t}\}| 
< k/2\bigr\}$, so a first order formula with parameters from $M$
\sn
\itemitem{ $(iii)$ }  {\rm Av\/}$_{\varphi(\bar y,\bar z)}
({\frak C},\bar{\bold b}^j)$ does not depend on $j$ if $\varphi(\bar y,
\bar z) \in \text{ {\rm avf\/}}(\bar{\bold b}^1) \cap \text{ {\rm
avf\/}}(\bar{\bold b}^2)$
\endroster}
\item "{$(c)$}"  {\rm avf\/}$(\bar{\bold b}^1) = \text{ {\rm avf\/}}
(\bar{\bold b}^2)$ so we can call it {\rm avf\/}$(p)$!
\sn
\item "{$(d)$}"  similarly $Cb_\varphi(\bar{\bold b}^1) = Cb_\varphi
(\bar{\bold b}^2)$ call it $Cb_\varphi(p)$ and $Cb(\bar{\bold b}_1) = 
Cb(\bar{\bold b}_2)$
call it $Cb(p)$, so $p \restriction \text{ {\rm avf\/}}(p)$ is definable
with parameters from $Cb(p)$
\sn
\item "{$(e)$}"  if $\varphi(\bar y,\bar z) \in \text{ {\rm avf\/}}(p)$, \ub{then}
$$
\text{{\rm Av\/}}_{\varphi(\bar y,\bar z)}({\frak C},\bar{\bold b}^1) = 
\text{ {\rm Av\/}}_{\varphi({\frak C},\bar y,\bar z)}
({\frak C},\bar{\bold b}^2).
$$
\endroster
\endproclaim
\bigskip

\demo{Proof}  Straight.
\mn
\ub{Clause (a)}:

There is no harm with increasing $I_j$, so \wilog, (if we supply
appropriate $\bar b_t$'s) $I_j$ has no last element.  We then prove by induction
on $\alpha$ (for all $A$).
\bn
\ub{Case 1}:  \ub{$\alpha = 0$}  Nothing to prove.
\bn
\ub{Case 2}:  \ub{$\alpha =1$}

Let $n < \omega$ be above the number of free variables in any formula in
$\Delta_1$ and let $t_0 < \ldots < t_{n-1} < t$ be in $I_j$.  Now
tp$_{\Delta_1}(\bar b_t,A \cup \dbcu_{\ell < n} \bar b_{t_\ell},{\frak C})$
is a finite set of formulas and let $\psi(\bar x,\bar c)$ be its conjunction.
So $\psi(\bar x,\bar c) \in \text{ tp}(\bar b_t,M \cup \{\bar b_s:s < t\},
{\frak C}) = \text{ {\rm Av\/}}(M \cup \{\bar b_s:s < t\},D_j)$ hence 
$\bold I = \{\bar b \in
{}^m M:{\frak C} \models \psi(\bar b,\bar c)\} \in D_j$.

So define ${\Cal F} = {\Cal F}_A$ by ${\Cal F}() = \bold I$, why is it as
required?  Clearly $\bar b' \in {\Cal F}()$ implies that $[t_{n-1} < t' \in
I_j \Rightarrow \bar b_{t'},\bar b'$ realizes the same $\Delta_1$-type over
$A \cup \{\bar b_t,\dotsc,\bar b_{t_{n-1}}\}]$.

But as $\langle \bar b_s:s \in I_j \rangle$ is an indiscernible sequence over
$M$, and $\bar b' \in M$, we can replace $t_0 < \ldots < t_{n-1}$, by any
$t'_0 < \ldots < t'_{n-1}$. \nl
But by the choice of $n$, for any $\bar b$ we have

$$
\align
\text{ tp}_{\Delta_1}(\bar b,A \cup \{\bar b_s:s \in I\}) = \cup
\{ &\text{tp}_{\Delta_1}(\bar b,A \cup \{\bar b_{s_0},\dotsc,b_{s_{n-1}}\}) \\
  &s_0 < \ldots < s_{n-1} \text{ are in } I_j\}
\endalign
$$
\mn
so we are done.
\bn
\ub{Case 3}: $\alpha > 1,\alpha < \omega$.

We define ${\Cal F}(\bar b_0,\dotsc,\bar b_{\ell -1})$ as
${\Cal F}_{A \cup\{\bar b_0,\dotsc,\bar b_{\ell -1}\}}$, and the checking is
easy.
\bn
\ub{Clause (b)}:

Subclause (i) holds by Claim \scite{np1.3}, Clause (d).

Subclause (ii) follows from subclause (i) as
\mr
\item "{$\boxtimes$}"  for $\bar c \in {}^{(\ell g(\bar z))}M$ we have:
$\varphi(x,\bar c) \in p$ iff 
$\varphi(\bar x,\bar c) \in \text{ {\rm Av\/}}(M,\bar{\bold b}^j)$ iff
\nl
$\varphi(\bar y,\bar c) \in \text{ av}_\varphi({\frak C},\langle \bar b_\ell:
\ell < k \rangle)$ iff ${\frak C} \models \vartheta[\bar c,b_0,\dotsc,
b_{k-1}]$ where \nl
$\vartheta(y,b_0,\dotsc,\bar b_{k-1}) = \dsize
\bigvee_{u \subseteq k,|u| \ge k/r} \,\, \dsize \bigvee_{\ell \in u}
\varphi(y,\bar b_\ell)$ \nl
[why?  the first ``iff" as Av$_\varphi({\frak C},\bar b^j)$ restricted to
$M$ is $p \restriction \varphi$, 
by an assumption of our claim, the second iff by clause (a) which we
have proved, the third iff by Definition \scite{np1.2}(2).]
\ermn
To Subclause (iii), we can choose for $j=1,2$, sequence $\bar b^j_\ell \in
{}^m M$ for $\ell < k$ as in clause (a), so for any $\bar c \in
{}^{(\ell g(z))} {\frak C}$ we have
\mr
\item "{$(*)$}"  $\varphi(\bar y,\bar c) \in \text{ Av}_{\varphi(\bar y,
\bar z)}({\frak C},\bar b^j)$ iff ${\frak C} \models \vartheta[\bar c,
\bar b^j_0,\dotsc,\bar b^j_{k-1}]$.
\ermn
So if the conclusion fails then for some $\bar c$ we have ${\frak C} \models
\Theta[\bar c,\bar b^1_0,\dotsc,b^1_{k-1}] =$ \nl
$\neg \Theta[\bar c,\bar b^2_0,
\dotsc,\bar b^2_{k-1}]$.  As $\bar b^j_\ell \in {}^m M,M \prec {\frak C}$ clearly
there is $\bar c' \subseteq M$ as above but by subclause (2)

$$
M \models \Theta[\bar c',b^j_0,\dotsc,b^j_{k-1}] \Leftrightarrow \varphi(\bar y,
\bar c') \in p.
$$ 
\mn
We get contradiction.
\bn
\ub{Clause (c)}:

This follows by \scite{np1.3} clause (d) and \scite{np1.5}, clause (a)
above.
\bn
\ub{Clause (d)}:

Similar to the proof of clause (c), using \scite{np1.4c}(2).
\bn
\ub{Clause (e)}:

This is (b)(iii) + (c).  \hfill$\square_{\scite{np1.5}}$
\enddemo
\bn
Of course
\demo{\stag{np1.5a} Observation}  For any 
$M \prec {\frak C}$ and $p \in S^m(M)$,avf$(p),Cb(p)$ are 
well defined as there are ultrafilters $D$ on ${}^mM$ 
such that Av$(M,D) = p$. 
\enddemo
\bn
\ub{\stag{np1.6} Observation}  [$T$ has NIP]  1) If 
$\bar{\bold b}$ is an infinite indiscernible set over $\emptyset$ 
(i.e. order immaterial), \ub{then} 
\mr
\item "{$(a)$}"  avf$(\bar{\bold b}) = L(T)$ (i.e. all formulas) 
\sn
\item "{$(b)$}"  if $\bar{\bold b}$ is an indiscernible sequence over $A$
then it is an indiscernible set over $A$;
\sn
\item "{$(b)^+$}"  if $\bar{\bold b}$ is a $\{\varphi(\bar x_0,\dotsc,
\bar x_k;\bar c)\}$-indiscernible sequence but is a 
$\Delta_\varphi$-indiscernible set (over $\emptyset$), \ub{then}
$\bar{\bold b}$ is a $\{\varphi(\bar x_0,\dotsc,\bar x_{k-1};
\bar c)\}$-indiscernible set when \nl
$\Delta_\varphi = \{\exists \bar x \dsize \bigwedge_{\ell < n}
\varphi(\bar x,\bar y_\ell)^{\text{if}\eta(\ell)}:\eta \in {}^n 2\}$ and $n$
is such that $\boxtimes^n_{\varphi(\bar x,\bar y)}$ from \scite{np1.0a}
fail. 
\ermn
2) If $p \in S^m(M)$ and $D_j,\bar{\bold b}^j$ for $j=1,2$ are as in
\scite{np1.5}, \ub{then} $\bar{\bold b}^1,\bar{\bold b}^2$ are nb-s of
distance 2.   
\bigskip

\demo{Proof} 1)(a) By \cite[Ch.II,4.13]{Sh:c} + the definitions. \nl

$(b)$   Easy by clause (a), see details in \scite{np4.4}(1). 
\nl

$(b)^+$  Similarly. \nl
2) By compactness and the proof of \scite{np1.5}.   \hfill$\square_{\scite{np1.6}}$
\enddemo
\bigskip

\demo{\stag{np1.7} Conclusion}  [$T$ has NIP]  In 
\scite{np1.5}, $\bar{\bold b}^1$ is an
indiscernible set over $\emptyset$ iff $\bar{\bold b}^1$ is an 
indiscernible set over $M$ iff $\bar{\bold b}^2$ is an indiscernible
set over $\emptyset$.
\enddemo
\bigskip

\definition{\stag{np1.8} Definition}  If $p \in S(M)$ and for some 
($\equiv$ all) $D^1,\bar{\bold b}^1$ as in \scite{np1.5} we have 
$\bar{\bold b}^1$ is an indiscernible set, \ub{then} we call 
$p$ a stable type.  Otherwise $p$ is called nonstable type.
\enddefinition
\bigskip

\demo{\stag{np1.8A} Conclusion}  [$T$ has NIP]  1) If 
$p \in S^m(M)$ is a stable type, \ub{then}
each $p \restriction \varphi$ is definable, in fact by parameters from
$Cb(p)$. \nl
2) The number of stable $p \in S^m(M)$ is $\le \|M\|^{|T|}$.
\enddemo
\bigskip

\demo{Proof}  1) By \scite{np1.5} (use clause (a) in \scite{np1.5}). \nl
2) Count the possible number of Definitions.  \hfill$\square_{\scite{np1.8A}}$
\enddemo
\bigskip

\remark{\stag{np1.8B} Remark}  
Note that $p \in S(M)$ may be definable but not stable, e.g. $M \prec
N$ are models of the theory of $(\Bbb R,<)$, and $a \in N
\backslash M$ is above all $b \in N$, then tp$(a,M,N)$ is definable but not
stable.
\endremark
\bn
\ub{\stag{np1.9} Observation}:  [$T$ is NIP and unstable]. \nl
1) There are $M \prec {\frak C}$ and nonstable $p \in S^1(M)$ [in ${\frak C}$
not ${\frak C}^{\text{eq}}$!]. \nl
2) There is an indiscernible sequence of \ub{elements} which is not an
indiscernible set of elements (over $\emptyset$!) \nl
3) If $\bar{\bold b} = \langle \bar b_t:t \in I \rangle$ is an indiscernible
sequence of $m$-tuples but not an indiscernible set, say for 
$\varphi(\bar x_0,\dotsc,\bar x_{k-1}),I$ a dense (linear order) for
simplicity \ub{then} for some $i < k-1$ and some $t_0,\dotsc,t_{k-1} \in I$
with no repetitions such that $t_i < t_{i+1},(t_i,t_{i+1})_I \cap 
\{t_0,\dotsc,t_{k-1}\} = \emptyset$ we have: the formula
$\psi(\bar b_{t_0},\dotsc,\bar b_{t_{i-1}},\bar x_i,x_{i+1},
\bar b_{t_{i+1}},\ldots)$ define a partial order on ${}^m{\frak C}$ by which
$\langle \bar b_t:t \in I,t_i <_I t <_I t_{i+1} \rangle$ is 
strictly increasing where

$$
\align
\psi(\bar x_{t_0},\dotsc,\bar x_{i-1},\bar x_i,\bar x_{i+1},\dotsc,
\bar x_{k-1}) = \forall \bar y[&\varphi(\bar x_{t_0},\dotsc,\bar x_{i-1},
\bar x_i,\bar y,x_{i+2},\dotsc,\bar x_{k-1}) \rightarrow \\
  &\varphi(\bar x_{t_0},\dotsc,\bar x_{i-1},\bar x_{i+2},\dotsc,
\bar x_{i+2},\dotsc,\bar x_{k-1})].
\endalign
$$
\bigskip

\demo{Proof}  1) As $T$ is unstable, for some $M \prec {\frak C}$ we have
$|S(M)| > \|M\|^{|T|}$ hence \scite{np1.8A}(2) some type $p \in S(M)$ is
nonstable. \nl
2) By part 1) and Definition \scite{np1.8}. \nl
3) Well known, see \cite[Ch.II,\S4]{Sh:c}  \hfill$\square_{\scite{np1.9}}$
\enddemo
\bn
Now \scite{np1.9}(3) applies to \scite{np1.9}(2) 
(where $\bar x_i$ is $x_i$) gives
\demo{\stag{np1.9A} Conclusion}  If $T$ is (NIP but) unstable, \ub{then} some
formula $\varphi(x,y;\bar c)$ define on ${\frak C}$ a quasi order with
infinite chains, (so $x,y$ singletons).
\enddemo
\bigskip

\remark{\stag{np1.10} Remark}  So if $T$ satisfies some version of $*$-stable
(see \cite[Ch.II]{Sh:300} or \cite{Sh:702}) \ub{then} 
$T$ is stable or $T$ has IP.
\endremark
\bn
So we may wonder \nl
\ub{Question}:  Does the ``has NIP" case in \scite{np1.10} is needed?  If
$T$ has IP does some $\varphi(x,y,\bar c)$ have the independence property?
\bn
Note that
\proclaim{\stag{np1.11} Claim}:  It $T$ is unstable, \ub{then} some formula
$\varphi(x,y,\bar c)$ has the order property (equivalently is unstable,
hence some $\varphi(x(y,\bar c))$ define a partial order with infinite claims
or has the independence property.
\endproclaim
\bigskip

\demo{Proof}  We know that some $\varphi(x,\bar y)$ is unstable so
choose a formula $\varphi(x,\bar y,\bar c)$ with the order property,
such that $\ell g(\bar y)$ is minimal.
So there is an indiscernible sequence $\langle \bar a_i \char 94 <b_i>:i <
\omega 4\rangle$ such that ${\frak C} \models \varphi[b_i,\bar a_j]$ iff
$j < i$.  Clearly $\langle b_i:i < \omega 4 \rangle$ is an indiscernible
sequence over $\bar c$, if it is not an indiscernible set, say not
$(\vartheta,k)$-indiscernible set, $\vartheta = \vartheta(x_0,\dotsc,
x_{k-1},\bar c)$, then for some permutation $\pi$ of $\{0,\dotsc,k-1\}$ and
$m$ the formula $\vartheta(x,y,\bar a_0,\dotsc,a_{m-1},a_{n-2},a_{2 \omega+1},
\dotsc,a_{2\omega+k,m-3},\bar c)$ linear orders $\langle a_{\omega +i}:i < \omega
\rangle$, hence has the order property.  So assume $\langle b_i:
i < \omega \rangle$ is an indiscernible set over $\bar c$, and let $a'_i$ be the first
element of the sequence $\bar a_i$.  If $\langle b_{2i+1}:i < \omega 4 \rangle$ is not an
indiscernible sequence over $\bar c \cup \{a'_{2 \omega}\}$ then we can find
a formula $\vartheta(x,y,\bar c'),\bar c' \subseteq \bar c \cup \{a_i:i <
\omega$ or $\omega 3 \le i\}$ such that $\models \vartheta[b_{2 \omega},
a_{\omega +2i+1}\bar c']$ for $i < \omega$ but $\models \neg \vartheta
[b_{2 \omega},a_{\omega 2 +2i+1}\bar c']$ for $i < \omega$ and we are 
done.  So assume
$\langle b_{2i+1}:i < \omega 4 \rangle$ is an indiscernible sequence over
$\bar c \cup \{a'_{2 \omega}\}$, hence all $\{a'_{2j}:j < \omega 4\}$
realizes the same type over $\{b_{2i+1}:i < \omega 4\} \cup \bar c$ hence for
$j < 2 \omega$ we can find $\bar a^*_{2j}$ realizing tp$(\bar a_{2j},
\{b_{2i+1}:i < \omega 4\} \cup \bar c,{\frak C})$ and the first element of
$\bar a^*_{2j}$ is $a'_0$.  This contradicts the choice of $\varphi(x,\bar y,
\bar c)$ as having the order property with $\ell g(\bar y)$ minimal as
we can ``move" $a'_0$ to $\bar c$.  \hfill$\square_{\scite{np1.11}}$
\enddemo
\bn
\centerline {$* \qquad * \qquad *$}
\bigskip

\remark{Remark}  Note that for indiscernible sets, the theorems on
dimension in \cite[III]{Sh:c} holds for theories $T$ with NIP, see \S3.
\endremark
\newpage

\head {\S2 Characteristics of types} \endhead  \resetall \sectno=2
\bigskip

We continue to speak on canonical bases and we deal with the
characteristics of tyeps and of indiscernible sets.  More elaborately,
for any indiscernible sequence $\bar{\bold b} = \langle \bar b_t:t \in
I \rangle,I$ an infinite linear order, we have a measure
ch$(\bar{\bold b}) = \langle CH_{\varphi(\bar x,\bar y_\varphi)}(\bar{\bold
b}):\varphi(\bar x,\bar y_\varphi) \in L(T) \rangle$ with $\bar
\lambda = \langle x_i:i < m \rangle,m = \ell g(\bar b_t)$ for $t \in
J$, where Ch$_{\varphi(\bar x,\bar y_\varphi)}(\bar{\bold b})$ measure how
badly $\varphi(\bar x,\bar y_\varphi)$ fail to be in avf$(\bar{\bold
b})$ (see Definition \scite{np2.5}), we can find such $\bar{\bold
b}$'s with maximal such CH$(\bar{\bold b})$ and wonder what can we say
about them.
\bigskip

\demo{\stag{np2.0} Hypothesis}  $T$ is NIP.
\enddemo
\bigskip

\definition{\stag{np2.1} Definition/Claim}  Let $\bar{\bold b} = \langle
\bar b_t:t \in I \rangle$ be an infinite indiscernible sequence, $k < \omega$.
\ub{Then}
\mr
\item "{$(a)$}"  (Claim) $\quad$ if $t_i \in I$ and 
$i < j \Rightarrow t_i <_I t_j$ for $i < j < \omega$ and \nl

$\qquad \qquad \quad \bar{\bold b}^k = 
\langle \bar b_{t_{ki}} \char 94 \bar b_{t_{ki+1}}
\char 94 \ldots \char 94 \bar b_{t_{ki+k-1}}:i < \omega \rangle$ then
{\roster
\itemitem{ $(\alpha)$ }  $Cb(\bar{\bold b}^1) \subseteq Cb(\bar{\bold b}^k)$,
\sn
\itemitem{ $(\beta)$ }  if $\varphi'(\bar x_1,\dotsc,\bar x_k;\bar y)
= \varphi(\bar x_\ell,\bar y)$ then: \nl

$\qquad \varphi'(\bar x_1,\dotsc,\bar x_k;\bar y) \in \text{
avf}(\bold b^k)$ \ub{iff} $\varphi(\bar x;\bar y) \in \text{ avf}(\bold{\bar b})$
\sn
\itemitem{ $(\gamma)$ }  if $\bar{\bold b}^{k,1},\bar{\bold b}^{k,2}$ are 
related like $\bar{\bold b}^k$ above to our $\bar{\bold b}$ then 
$Cb(\bar{\bold b}^{k,1}) = Cb(\bar{\bold b}^{k,2})$
\sn
\itemitem{ $(\delta)$ }  if $\varphi(\bar x,\bar y) \in \text{ daf}(\bold{\bar b}),
\varphi' = \varphi'(\bar x_1,\dotsc,\bar x_k,\bar y) = \varphi(\bar
x_\ell,\bar y) \and \neg \varphi(\bar x_m,\bar y)$ or $= \varphi(\bar
x_\ell,\bar y) \equiv \neg \varphi(\bar x_m,\bar y)$ then $\varphi'
\in \text{ avf}(\bar b)^k$
\endroster}
\sn
\item "{$(b)$}"   (Definition) $\quad$ let $Cb^k(\bar{\bold b}) =
Cb(\bar{\bold b}^k)$,Av$^k(\bar{\bold b},{\frak C}) = \text{ Av}_{\text{avf}}
(\bar{\bold b}^k,{\frak C})$ for any $\bar{\bold b}^k$ \nl

$\qquad \qquad \,\,\,$ as above
\sn
\item "{$(c)$}"  (Definition) $\quad Cb^\omega(\bar{\bold b}) = \cup
\{Cb^k(\bar{\bold b}):k < \omega\}$
\sn
\item "{$(d)$}"  (Fact) $\quad$ if $I_1,I_2$ are infinite subsets of $J$ and
$\bar{\bold b} = \langle \bar{\bold b}_t:t \in J \rangle$ an \nl

$\qquad \qquad \quad$ indiscernible
sequence (recall $J$ linear order) \ub{then} \nl

$\qquad \qquad \quad Cb^\omega(\bar{\bold b} \restriction I_1) = 
Cb^\omega(\bar{\bold b} \restriction I_2)$.
\sn
\item "{$(e)$}"  (Fact) $\quad$ If the infinite indiscernible
sequences $\bar{\bold b}^1,\bar{\bold b}^2$ are nb-s, \ub{then}
Cb$^\alpha(\bar{\bold b}^1) = \text{ Cb}^\alpha(\bar{\bold b}^2)$ for
$\alpha \le \omega$.
\endroster
\enddefinition
\bigskip

\demo{Proof}  Easy.
\enddemo
\bigskip

\definition{\stag{np2.2} Definition}  For $\alpha \le \omega$.
We say $p \in S^m(A)$ does not $\alpha$-fork over
$B \subseteq A$, \ub{if} for some model $M \supseteq A$ and $q \in S^m(M)$
extending $p$ we have $Cb^\alpha(q) \subseteq acl_{{\frak C}^{\text{eq}}}(B)$.
Similarly we say that $C/B$ does not $\alpha$-fork over $A \subseteq B$ if 
$\bar c \subseteq C \Rightarrow \text{ tp}(C,B)$ does not $\alpha$-fork 
over it.
\enddefinition
\bigskip

\proclaim{\stag{np2.3} Claim}  1) In \scite{np2.2}: ``for some $M \supseteq A$"
can be replaced by ``for every $M \supseteq A$".
\endproclaim
\bigskip

\demo{Proof}  Easy.
\enddemo
\bigskip

\remark{\stag{np2.4} Remark}:  Assume that $T$ is a simple theory, 
$\bar{\bold b} = \langle \bar b_t:t \in I \rangle$ is an infinite 
indiscernible sequence.  \ub{Then} 
we cannot find $\langle \bar a_n:n < \omega \rangle$ indiscernible, 
$\langle \varphi(\bar x,\bar a_n):n < \omega \rangle$ pairwise contradictory 
(or just $m$-contradiction for some $m$) and

$$
\dsize \bigwedge_{n < \omega} (\exists^\infty t \in I)(\varphi(\bar b_t,
\bar a_n).
$$
\endremark
\bigskip

\demo{Proof}   As we can repeat and get the tree property.
More fully, for any cardinals $\mu > \kappa$ we consider $J =
{}^\kappa \mu$ as a linearly ordered set,
ordered lexicographically and for $\rho \in {}^{\kappa >} \mu$ let
$J_\rho = \{\nu \in J:\rho \triangleleft \nu\}$; \wilog \, $I$ is
countable and $h:I \rightarrow J$ is order preserving.  We can find $\bar a_\eta
\in {\frak C}$ for $\eta \in J$ such that $\langle c_\eta:\eta \in J
\rangle$ is an indiscernible sequence such that $t \in I \Rightarrow
c_{h(t)} = b_t$.  By compactness, for each $\alpha < \kappa$ we can
find $\langle a_\rho;\rho \in {}^\alpha \mu \rangle$ such that:
\mr
\item "{$(\alpha)$}"  $\langle \varphi(\bar x,\bar a_\rho):\rho \in
{}^\alpha \mu \rangle$ are pairwise contradictory (or just any $m$ of
them)
\sn
\item "{$(\beta)$}"  $\eta \in J_\rho,\rho \in {}^\alpha \mu
\Rightarrow {\frak C} \models \varphi[c_\eta,a_\rho]$.
\ermn
Now $\langle \varphi(\bar x,\bar a_\rho):\rho \in {}^{\kappa >} \mu
\rangle$ exemplified the tree property. \hfill$\square_{\scite{np2.4}}$
\enddemo
\bn
We have looked at indiscernible sequences which are stable.  We now look
after indiscernible sequences which are in the other extreme.
\definition{\stag{np2.5} Definition}  For $\bar{\bold b} = \langle \bar b_t:
t \in I \rangle$ an indiscernible sequence, we define its character

$$
Ch(\bar{\bold b}) = \langle Ch_{\varphi(\bar y,\bar z)}(\bar{\bold b}):
\varphi(\bar y,\bar z) \in L(T) \rangle
$$
\mn
where

$$
\align
Ch_{\varphi(\bar y,\bar z)}(\bar{\bold b}) = \text{ Max}\{n:&\text{for some }
\bar c,\langle \text{ TV}(\varphi(\bar b_t,\bar c):t \in I \rangle \\
  &\text{change sign } n \text{ times (i.e. } I \text{ divided to } n+1
\text{ intervals})\}.
\endalign
$$
\mn
2) For $p \in S^m(A)$, let
\mr
\item "{$(a)$}"  CH$(p) = \{Ch(\bold{\bar b}):\bar{\bold b} \text{ is
an infinite indiscernible sequence such that every } \bar b_t \text{
realizes } p)\}$
\sn
\item "{$(b)$}"  for a formula $\varphi = \varphi(\bar x_0,\dotsc,\bar
x_{k-1})$ let CH$(p,\varphi(\bar x_0,\dotsc,x_k) =
\{\text{Ch}(\bar{\bold b}):\bar{\bold b} = \langle \bar b_t:t \in I
\rangle$ is an infinite indiscernible set such that $t_0 <_I t_1 <_I
\ldots <_I t_{k-1} \Rightarrow {\frak C} \models \varphi[\bar
b_{t_0},\dotsc,\bar b_{t_{k-1}}]$
\sn
\item "{$(c)$}"  CH$^{\text{max}}(p) = 
\{\bar n \in CH(p):\text{there is no bigger such } \bar n' \in
CH(p)\}$, when ``$\bar n'$ is bigger than $\bar n$" mean $(\forall \varphi)
(n_\varphi \le n'_\varphi) \and (\exists \varphi)(n_\varphi <
n'_\varphi)$
\sn
\item "{$(d)$}"  CH$^{\text{min}}(p,\varphi(\bar x_0,\dotsc,x_{k-1})) =
\{\bar n \in \text{ CH}(p,\varphi(\bar x_0,\dotsc,\bar x_{k-1}))$:
there is no smaller $\bar n' \in \text{ CH}(p,\varphi(\bar
x_0,\dotsc,\bar x_{k-1}))\}$.
\ermn
Note: for the trivial $\varphi$, CH$(p,\varphi) = \text{ CH}(p)$ hence
CH$^{\text{max}}(p,\varphi) = \text{ CH}^{\text{max}}(p)$.
\enddefinition
\bigskip

\proclaim{\stag{np2.6} Claim}  Let $p \in S^m(A)$ be non-algebraic, $\bar x =
\langle x_\ell:\ell < n \rangle$. \nl
1) If $\bar n = \langle n_\varphi:\varphi = \varphi(\bar x,\bar y) \rangle
\in CH(p)$, \ub{then} there is $\bar n' \in CH^{\text{ max}}(p)$ such that
$\bar n \le \bar n'$. \nl
2) CH$^{\text{max}}(p)$ is non-empty. \nl
3) If CH$(p,\varphi) \ne \emptyset$ then CH${\text{min}}(p,\varphi)
\ne \emptyset$ and CH${\text{max}}(p,\varphi) \ne \emptyset$. 
\endproclaim
\bigskip

\demo{Proof}  Let $R,<$ be an $n$-place and $2n$-place predicate not in
$\tau_T$ and let

$$
\align
\Gamma_p = Th({\frak C}_T,c)_{c \in A} &\cup \{(\forall \bar x)[R(\bar
x) \rightarrow \Theta(\bar x),\bar c)]:\Theta(\bar x,\bar c) \in p\}
\\
  &\cup \{(\exists \bar x_0,\dotsc,
\bar x_{n-1})(\dsize \bigwedge_{\ell < k} R(\bar x_\ell) \and 
\dsize \bigwedge_{\ell < m} \bar x_\ell \ne \bar x_m):n < \omega\} \\
  &\cup \{\bar x < \bar y \rightarrow R(\bar x) \wedge R(y)\} \\
  &\cup \{``\text{and } \le \text{ linearly ordered } \{\bar x:R(\bar x)\}"\} \\
  &\cup \{(\forall \bar x_1),\dotsc,
(\forall \bar x_m)(\forall \bar y_1) \ldots (\forall \bar y_m)(\bar x_1 <
\bar x_2 < \ldots < \bar x_m \\
  &\and \bar y_1 < \ldots < \bar y_m 
\rightarrow \psi(\bar x_1,\dotsc,\bar x_m,\bar c) 
\equiv \psi(\bar y_1,\dotsc,\bar y_m,\bar c): \\
  &\bar c \subseteq A \text{ and } \psi \in L(\tau_T) \text{ and }
m < \omega\} 
\endalign
$$
\mn
(with $\bar x_i = \langle x_{i,\ell}:\ell < m \rangle$).  If
$\lambda = |T| + \aleph_1$ we may omit it.  For $\bar n = \langle 
n_{\varphi(\bar x,\bar y)}:\varphi(\bar x,\bar y) \in L(\tau_T) \rangle$ and 
$\bar \varphi = \langle \varphi_i(\bar x,\bar y_i):i < |T| \rangle$ listing
those formulas let

$$
\align
\Gamma_{\bar n,\bar \varphi} = &\{\vartheta_{n_i,\varphi_i}:i < |T|\}
\text{ where} \\
  &\vartheta_{n,\varphi(\bar x,\bar y)} = (\exists \bar y)
(\exists \bar x_0,\dotsc,\exists \bar x_n)[\bar x_0 < \bar x_1 < \ldots
< \bar x_n \and \\
  &\hskip80pt \dsize \bigwedge_{\ell < n} (\varphi(\bar x_\ell,\bar y))
\equiv \neg \varphi(\bar x_{\ell +1},\bar y)].
\endalign
$$
\mn
Now easily
\mr
\item "{$(a)$}"  $\Gamma_p$ is a consistent type (using $p$ being
non algebraic and Ramsey theorem)
\sn
\item "{$(b)$}"  $\Gamma_p \cup \Gamma_{\bar n,\bar \varphi}$ is consistent
iff $\bar n \in \text{ CH}(p)$
\sn
\item "{$(c)$}"  if $\bar n \le \bar n'$ then $\Gamma_{\bar n,\bar \varphi}
\subseteq \Gamma_{\bar n',\bar \varphi}$
\sn
\item "{$(d)$}"  if $J$ is a directed order 
$\bar n_t = \langle n_{t,\varphi(\bar x,\bar y_v)}:\varphi(\bar x,\bar
y_\varphi) \rangle \in \text{ CH}(p,\varphi)$ increases with $t \in J$
and $\bar n^* = \langle n^*_{\varphi(\bar x,\bar
y_\varphi)}:\varphi(x,\bar y_\varphi) \rangle,n^* = \text{
max}\{n_{t,\varphi(\bar x,\bar y_\varphi)}:t \in J\}$, \ub{then} $n^*
\in \text{ CH}(\bar p^*,\varphi)$
\sn
\item "{$(e)$}"  like (d) investing the order.
\ermn
Together we can deduce the desired conclusions.
\hfill$\square_{\scite{np2.6}}$
\enddemo
\bn
\ub{\stag{np2.7} Question}:   For $p \in S(A)$ (or just $p \in S(M))$ does
indiscernible sequences $\bar{\bold b}$ of elements realizing $p$ such that
Ch$(\bar{\bold b}) \in \text{ CH}^{\text{max}}(p)$, CH$(\bar b) \in
\text{ CH}^{\text{min}}(p,\varphi)$ play a special role?
\newpage

\head {\S3 Shrinking indiscernibles} \endhead  \resetall \sectno=3
\bn
The case of indiscernible sets is easier so we ignore it.
\bigskip

\proclaim{\stag{3.4} Claim}  If $\bar{\bold b}= \langle \bar{\bold b}_t:t \in
I \rangle$ is an indiscernible sequence over $A,\bar c \in {\frak C}$ (so
finite), \ub{then} 
\mr
\item "{$(a)$}"  there are $J \subseteq I,J^* \subseteq J,|J^*| 
\le |T|$
such that
{\roster
\itemitem{ $(*)$ }  if $n < \omega,\bar s,\bar t \in {}^n I,\bar s \sim_{J^*} \bar t$
(i.e. $\bar s,\bar t$ realize the same quantifier free type in the
linear order $J$) \ub{then} $\bar a_{\bar s} = \langle
\bar a_{s_\ell}:\ell < n \rangle,\bar a_{\bar t} = \langle
a_{t_\ell}:\ell < n \rangle$ realize in ${\frak C}$ the
same type over $A \cup \bar c$
\endroster}
\item "{$(b)$}"  if we fix $n$ and deal with $\varphi$-types we can
demand $|J^*| < k_{\varphi,n} < \omega$
\sn
\item "{$(c)$}"  if in adition $\bar{\bold b}$ is an indiscernible
set, \ub{then} in $(*)$ of clause (a) we cn weaken $\bar s \sim_{J^*}
\bar t$ to $(\forall \ell,k)[(s_\ell < s_k \equiv t_\ell < t_k) \and
s_\ell \in J^* \equiv t_\ell \in J^* \rightarrow s_\ell = t_\ell]$.
\endroster
\endproclaim
\bigskip

\demo{Proof}

$(a) \quad$ by $(b)$ \nl

$(b) \quad$ follows by Claim \scite{npx.2} below

$(c) \quad$ similarly.  \hfill$\square_{\scite{3.4}}$
\enddemo
\bn
\centerline {$* \qquad * \qquad *$}
\bigskip

\definition{\stag{npx.1} Definition}  1) For a linear order $I,m^* \le
\omega,n \le \omega,\alpha_\ell$ an ordinal, a model $M$ and a set 
$A \subseteq M$, we say that $\bold{\bar a} = 
\langle a_{u,\alpha,\ell}:\ell < n,u \in [I]^\ell,
\alpha < \alpha_{|u|} \rangle$ is $(\Delta^*,m^*)$-indiscernible over $A$
of the $\langle \alpha_\ell:\ell < n \rangle$-kind if the following holds:
\mr
\item "{$(*)$}"  if $m \le m^*,I \models t_0 < \dots < t_{m-1},I \models
s_0 < \ldots < s_{m-1}$ for $v \subseteq m$ we let $u_v = \{t_\ell:\ell \in
v\},w_v = \{s_\ell:\ell \in v\}$ then $\langle a_{u_v,\alpha,\ell}$: for 
$\ell \le m,\ell < n,v \in [m]^\ell,\alpha < \alpha_\ell \rangle$ and $\langle
a_{w_v,\alpha,\ell}:\ell \le m,v \in [m]^\ell,\alpha < \alpha_\ell \rangle$
realizes the same $(\Delta,m^*)$-type over $A$ in $M$.
\ermn
2) If we omit $\Delta$ we mean all first order formulas, if we omit $m^*$
we man $\omega$.  Also in $a_{u,\alpha,\ell}$ we may omit $\ell$ (it is
$|u|$).  Of course nothing changed if we allow $a_{u,\alpha,\ell}$ to be a
finite sequence (with length depending on $(\alpha,\ell)$ only. \nl
3) We add ``and over $J$" where $J \subseteq I$ if in 
$(*)$ we demand $(\forall x \in J)
\dsize \bigwedge_\ell(x < t_\ell = x < s_\ell \and x = t_\ell 
\equiv x = s_\ell \and t_\ell < x \equiv x_\ell < x)$.  We say ``almost over
$J$" if we add $J \cap \{t_\ell:\ell < n\} = \emptyset$.
\enddefinition
\bigskip

\proclaim{\stag{npx.2} Claim}  [$T$ has NIP]  1) Assume
\mr
\item "{$(a)$}"  $\Delta$ is a finite set of formulas, $m^* < \omega$
\sn
\item "{$(b)$}"  $M$ a model of $T,A \subseteq M$
\sn
\item "{$(c)$}"  $\bold a = \langle a_{u,k,\ell}:\ell < n,\alpha < k_\ell,
u \in [I]^\ell \rangle$ is indiscernible over $A$
\sn
\item "{$(d)$}"  $\bar d \in {}^{\omega >}M$.
\ermn
\ub{Then} there is a finite subset $J$ of $I$ such that $\langle
a_{u,k,\ell}:\ell < n,k < k_\ell,u \in [I^\ell] \rangle$ is 
$\Delta$-indiscernible over $A \cup \bar d$ almost over $J$. \nl
2) Moreover, there is a bound on $|J|$ which depend just on $\Delta,
\langle k_\ell:\ell < n \rangle$ (and $T$), and so it is enough that $\bold a$
is $\Delta_1$-indiscernible for appropriate finite $\Delta_1$.
\endproclaim
\bigskip

\demo{Proof}  1) Straightforward.  If this fails, try to choose by induction
on $i < \omega$, \nl
$\langle t^i_\ell:\ell < m_i \rangle,
\langle s^i_\ell:\ell < m^*_i \rangle$ such that:
\mr
\widestnumber\item{$(iii)$}
\item "{$(i)$}"  $m_i \le m^*$
\sn
\item "{$(ii)$}"  $t^i_0 < t^i_1 < \ldots < t^i_{m_i-1},s^i_0 < s^i_1 <
\ldots < s^i_{m_i-1}$
\sn
\item "{$(iii)$}"  $t^i_m,s^i_m \notin J_i = \{t^j_\ell,s^j_\ell;j < i,\ell
< m_j\}$
\sn
\item "{$(iv)$}"  $\langle t^i_\ell:\ell < m_i \rangle,\langle s^i_\ell:
\ell < m_i \rangle$ exemplify that $J_i$ is not as required.
\ermn
Let $\bar b^0_i = \langle a_{u,k,\ell}:\ell < n,k < k_\ell,u \in
[\{t_0,\dotsc,t_{m_i-1}\}]^\ell \rangle$ and
$b^1_i = \langle a_{u,k,\ell}:\ell < n,k < k_\ell,u \in [\{s_0,\dotsc,
s_{m_i-1}\}]^\ell \rangle$.
\sn
So clearly
\mr
\item "{$(*)_1$}"  the $\Delta$-types of $\bar d \char 94 \bar b^0_i,
\bar d \char 94 b^{-1}_i$ over $A$ are different \nl
[why?  by their choice]
\sn
\item "{$(*)_2$}"  if $i(*) < \omega,\eta \in {}^{i(*)}2$, then the
types of $\bar b^0_0 \char 94 \bar b^0_1 \char 94 \ldots \bar b^0_{i(*)-1}$
and \nl
$\bar b^{\eta(0)}_0 \char 94 \bar b^{\eta(1)}_1 \char 94 \ldots
\char 94 b^{\eta(i(*)-1)}_{i(*)-1}$ over $A$ are equal \nl
[why?  by the indiscernibility].
\ermn
So we are easily done.  \hfill$\square_{\scite{npx.2}}$
\enddemo
\bigskip

\proclaim{\stag{npx.3} Claim}  [$T$ has NIP]  Assume 
$\bar{\bold a}^\ell = \langle
\bar a^\ell_t:t \in I_\ell \rangle$ is an indiscernible sequence of
cofinality $\kappa > |T|$ for $\ell = 1,2$.  \ub{Then} we can find
$s^\ell_i \in I_\ell$ for $\ell = 1,2,i < \kappa$ such that
$\langle a^1_{s^1_i} \char 94 \bar a^2_{s^2_i}:i < \kappa \rangle$ is
an indiscernible sequence.
\endproclaim
\bigskip

\demo{Proof}  Easy by repeated use of \scite{3.4}.
\enddemo
\bn
\ub{\stag{npx.4} Conclusion}  1) Assume
\mr
\item "{$(*)$}"  $\langle \bar b_t:t \in I \rangle$ is an indiscernible
sequence over $A$.
\ermn
For every $\bar c \in {}^{\omega >}{\frak C}$ there is $J \subseteq I$ of
cardinality $\le |T|$ and $\langle J_\varphi:\varphi \in L_{\tau(T)} \rangle,
J_\varphi$ a finite subset of $J$ such that:
\mr
\item "{$(*)_1$}"  for every $\bar a \in {}^{\ell g(\bar y)} A$ and
$\varphi = \varphi(\bar x,\bar y)$ there are $n \le 
n_{\varphi(\bar x,\bar y)}$ and $t_1 < \ldots < t_n$ from $J_\varphi$ such 
that if $r,s \in I \backslash \{t_1,\dotsc,t_n\}$ and $m \in [1,m] \Rightarrow
s <_I t_m \equiv r <_I t_n$ then $\models \varphi[\bar b_s,\bar a] \equiv
[\bar b_r,\bar a]$
\sn
\item "{$(*)_2$}" for every $k < \omega,\bar a \in {}^{\ell g(\bar y)} A$
and $\varphi = \varphi(\bar x_1,\dotsc,\bar x_k,\bar y)$ there are $n \le
n_\varphi$ and $t_1 < \ldots < t_n$ from $J_\varphi$ we have if $s_1 <_I
\ldots <_I s_k$ and $r_1 <_I \ldots <_I r_k$ are from $J$ and $m \in [1,n]
\and \ell \in [1,k] \Rightarrow (s_\ell <_I t_m \equiv r_\ell <_I t_m) \and
t_m <_I s_\ell \equiv t_m <_I r_\ell)$ then $\models \varphi[\bar b_{s_1},
\dotsc,\bar b_{s_k},\bar a] \equiv \varphi[\bar b_{r_1},\dotsc,\bar b_{r_k},
\bar a]$.
\ermn
2) Assume
\mr
\item "{$(*)_3$}"  $\langle \bar b_{u,\alpha,\ell}:\ell < n,u \in [I]^n,
\alpha < \alpha_\ell \rangle$ is indiscernible over $A$ and $\alpha_\ell <
\omega$ for $\ell < n$ (and $n < \omega$).  For every $\bar c$ there are
$J \subseteq J,|J| \le |T|$ 
and finite $J_\varphi \subseteq J$ for $\varphi \in
L_{\tau(T)}$ such that the parallel of $(*)_1,(*)_2$ hold. 
\endroster
\bigskip

\demo{Proof}  1) Clearly $(*)_1$ is a case of \scite{npx.2}, if we apply it
to $({\frak C},a)_{a \in A}$.  Similarly for $(*)_2$ apply it enough times
noting: if $\langle J_\ell:\ell \le k+1 \rangle$ is an increasing sequence of
subsets of $I$ and $t,\dotsc,t_k \in I$ then for some $\ell < k+1,\{t_1,
\dotsc,t_k\} \cap J_{\ell +1} \subseteq J_\ell$. \nl
2) Similar.
\enddemo
\bn
\ub{Question}:  Can we find $\bar b_\alpha \in {}^{|T|}{\frak C}$ such that
$\langle \bar b_\alpha:\alpha < \lambda \rangle$ is an indiscernible sequence
$\alpha \ne \beta \Rightarrow \bar b_\alpha \ne \bar b_\beta$ and for
$\alpha < \beta < \gamma$ we have tp$(\bar b_\gamma,\bar b_\beta) \models
\text{ tp}(\bar b_\gamma,\bar b_\alpha)$?
\bn
\ub{Question}:  If $< (= \varphi(x,y,\bar c))$ is a partial order with
infinite increasing sequences, we may consider $\kappa$-directed subsets,
$\kappa = \text{ cf}(\kappa > |T|)$, they define a Dedekind cut.

What about orthogonality of those?
\newpage

\head {\S4 Perpendicular endless indiscernible sequences} \endhead  \resetall \sectno=4
\bigskip

Dimension and orthogonality play important role in \cite{Sh:c}, see in
particular Ch.V.  Now, as our prototype is the theory Th$(\Bbb Q,<)$,
it is natural to look at cofinality, this is dual-cf$(\bar{\bold
b},A)$, measuring the cofinality of approaching $\bar{\bold b}$ from
above (here $\bar{\bold b}$ is always indiscernible sequences with no
lat member).  So a relative of orthogonality which we all
perpencidularity suggest itelf as relevant.  It is defined in
\scite{np4.1}, as well as equivalence and dual-cf.  Now
perpendicularity is closely related to being mutual indiscernibility
(see \scite{np4.2}(1), \scite{np4.3}(2), hence if $T$ is unstable,
then there are lost if pairwise perpendicular indiscernible sequences:
cf$\langle \bar a_\alpha:\alpha < \lambda \rangle$ is an indiscernible
sequence, not set and $\bar{\bold b}^\alpha = \langle \bar a_{\omega
\alpha + n}:n < \omega \rangle$ for $\alpha < \lambda$ then
$\{\bar{\bold b}^\alpha:\alpha < \lambda\}$ are pairwise
perpendicular.  In this section we present basic properties of
perpendicularity.  In particular, it is preserved by equivalence
(\scite{np4.3}(5)).  For perpendicular sequences, we can more easily
restrict them to get mutually indiscernible sets than in \S3.

For indiscernible sets this essentially becomes orthogonality.

The case of looking at more than two indiscernible sequences reduced
to looking at all pairs (\scite{np4.5}(2), \scite{np4.7}(2)).  Also,
as in \cite[V]{Sh:c}, if $\bar{\bold b}$ is not perpendicular to
$\bar{\bold a}^\zeta$ for $\zeta < \zeta^*$ and the $\bar{\bold
a}^\zeta$-s are pairwise perpendicular then $\zeta^* < |T|^+$ (see
\scite{np4,8a}).  

Lastly, we recall the density of quite ``types not splitting over
small sets" (for theories with NIP), hence the existence of a ``quite
constructible" model over any $A$.
\bigskip

\definition{\stag{np4.0a} Definition}  1) We say the infinite sequences
$\bar{\bold b}^1,\bar{\bold b}^2$ are mutually indiscernible if
$\bar{\bold b}^\ell$ is indiscernible over $\cup\{ \bar b^{3 -
\ell}_t:t \in \text{ Dom}(\bar{\bold b}^{3 - \ell})\}$ for $\ell
=1,2$.  Similarly over $A$. \nl
2) We say that the family $\{\bar{\bold b}^\zeta:\zeta < \zeta^*\}$ of
sequences is mutually indiscernible over $A$, \ub{if} for $\zeta <
\zeta^*,\bar{\bold b}^\zeta$ is indiscernible over $\cup\{ \bar
b^\varepsilon_t:\varepsilon \ne \zeta,\varepsilon < \zeta^*,t \in
\text{ Dom}(\bar{\bold b}^\varepsilon)\} \cup A$.
\enddefinition
\bigskip

\demo{\stag{np4.0} Hypothesis}  $T$ has NIP.
\enddemo
\bigskip

\definition{\stag{np4.1} Definition}  Let $\bar{\bold a}^\ell = \langle
a^\ell_t:t \in I_\ell \rangle$ be an indiscernible sequence which are
endless (i.e. $I_\ell$ having no last element) for $\ell = 1,2$. \nl
1) We say that $\bar{\bold a}^1,\bar{\bold a}^2$ are perpendicular when
\mr
\item "{$(*)$}"  \ub{if} $\bar b^\ell_n$ realizes Av$(\{\bar b^k_m:m < n
\and k \in \{1,2\} \vee m = n \and k < \ell\} \cup \bar{\bold a}^1 \cup
\bar{\bold a}^2,\bar{\bold a}^\ell)$ for $\ell =1,2$ then $\bar{\bold b}^1,
\bar{\bold b}^2$ are mutually indiscernible (see below) where
$\bar{\bold b}^\ell = \langle \bar b^\ell_n:n < \omega \rangle$.
\ermn
2) We say $\bar{\bold a}^1,\bar{\bold a}^2$ are equivalent if for every
$A \subseteq {\frak C}$ we have Av$(A,\bar{\bold a}^1) = \text{ Av}
(A,\bar{\bold a}^2)$. \nl
3) If $\bar{\bold a}^1 \subseteq A$ we let dual-cf$(\bar{\bold a}^1,A) =
\text{ Min}\{|B|:B \subseteq A$ and no $\bar c \in {}^{\omega >}A$ realizes
Av$(B,\bar{\bold a}^1)\}$.
\enddefinition
\bigskip

\proclaim{\stag{np4.2} Claim}  1) If $\bar{\bold a}^1,\bar{\bold a}^2$ are
endless mutually indiscernible, \ub{then} they are \nl
perpendicular. \nl
2) ``Mutually indiscernible" and ``perpendicular" are symmetric relations. \nl
3) On the family of endless indiscernible sequences, being equivalent is 
an equivalence relation.
\endproclaim
\bigskip

\demo{Proof}  Easy.
\enddemo
\bigskip

\proclaim{\stag{np4.3} Claim}  1) If $\bar{\bold a}^\ell = \langle a^\ell_t:
t \in I_\ell \rangle$ is an indiscernible sequence for $\ell =1,2$ and $|T| <
\text{ cf}(I_1),|I_1| < \text{ cf}(I_2)$, \ub{then} for some end segments
$J_1,J_2$ of $I_1,I_2$ respectively, $\bold{\bar a}^1 \restriction J_1,
\bold{\bar a}^2
\restriction J_2$ are mutually indiscernible. \nl
2) If $\bar{\bold a} = \langle a^\ell_t:t \in I_\ell \rangle$ is an
indiscernible sequence for $\ell=1,2$ and cf$(I_1)$, cf$(I_2)$ are 
infinite and
distinct \ub{then} $\bar{\bold a}^1,\bar{\bold a}^2$ are perpendicular. \nl
3) If $\bar{\bold a}^\ell = \langle \bar a^\ell_t:t \in I_\ell \rangle$
is an endless indiscernible sequence for $\ell = 1,2$, $\delta$
is limit ordinal and $\bar b^\ell_\alpha$ realizes Av$(\{\bar b^k_\beta:
\beta < \alpha \and k \in \{1,2\}$ or $\beta = \alpha \and k < \ell\} \cup
\bar{\bold a}^1 \cup \bar{\bold a}^2,\bar{\bold a}^\ell)$ and
$\bar{\bold b}^\ell = \langle \bar b^\ell_\alpha:\alpha < \delta \rangle$ for
$\ell = 1,2$ \ub{then}: $\bar{\bold a}^1,\bar{\bold a}^2$ are perpendicular
iff $\bar{\bold b}^1,\bar{\bold b}^2$ are perpendicular. \nl
4) If $\bar{\bold a}^\ell = \langle a^\ell_t:t \in I_\ell \rangle$ is an
endless indiscernible sequence and $J_\ell \subseteq I_\ell$ is unbounded
for $\ell = 1,2$, \ub{then} $\bar{\bold a}^1,\bar{\bold a}^2$ are perpendicular iff
$\bar{\bold a}^1 \restriction J_1,\bar{\bold a}^2 \restriction J_2$ are
perpendicular. \nl
5) If $\bar{\bold a}^\ell = \langle a_t:t \in I^\ell \rangle$ are an endless
indiscernible sequence for $\ell=1,2,3,4$ and $\bar{\bold a}^1,
\bar{\bold a}^3$ are equivalent and $\bar{\bold a}^2,\bar{\bold a}^4$ are
equivalent, \ub{then} $\bar{\bold a}^1,\bar{\bold a}^2$ are perpendicular
iff $\bar{\bold a}^3,\bar{\bold a}^4$ are perpendicular.
\endproclaim
\bigskip

\demo{Proof}  Straight.
\enddemo
\bigskip

\remark{Remark}  1) Replace $\bar{\bold a}$ by a sequence of concatanation
of $n$-tuples from it (as in \scite{np2.1}) preserve relevant properties.
\nl
2) In \scite{np4.3}(1), can we weaken $|I_1| < \text{ cf}|I_2|$ to
cf$|I_1| \ne \text{ cf}|I_2|$?
\endremark
\bigskip

\proclaim{\stag{np4.4} Claim}  Assume $\bar{\bold a}^\ell = \langle
\bar a^\ell_t:t \in I_\ell \rangle$ is an endless indiscernible sequence for
$\ell = 1,2$. \nl
1) If $\bar{\bold a}^1$ is an indiscernible sequence over $A$, \ub{then}:
$\bar{\bold a}^1$ is an indiscernible set over $A$ \ub{iff} $\bar{\bold a}^1$
is an indiscernible set over $\emptyset$. \nl
2) $\bar{\bold a}^1$ is nonstable in ${\frak C}$ \ub{iff} $\bar{\bold a}^1$ is
nonstable in $({\frak C},c)_{c \in A}$. \nl
3) If $\bar{\bold a}^1,\bar{\bold a}^2$ are equivalent, \ub{then}
$\bar{\bold a}^1$ is nonstable iff $\bar{\bold a}^1$ is nonstable. \nl
4) If $J_\ell \subseteq I_\ell$ is infinite and $\bar{\bold a}^1,
\bar{\bold a}^2$ are mutually indiscernible then $\bar{\bold a}^1
\restriction J_1,\bar{\bold a}^2 \restriction J_2$ are mutually indiscernible
over $\dbcu^2_{\ell = 1}(\bar{\bold a}^\ell \restriction (I_\ell \backslash
J_\ell))$.
\endproclaim
\bigskip

\demo{Proof}  1) The ``only if" is trivial.  For the other direction if
$\bar{\bold a}^1$ is an indiscernible set over $\emptyset$ but not over
$A$, we easily get the independence property.  I.e. assume that 
$t_0 < \ldots <
t_{n-1}$ in $I_1$ and $\pi$ a permutation of $\{0,\dotsc,n-1\},\varphi
(\bar b,\bar a_{t_0},\dotsc,\bar a_{t_{n-1}}) \and \neg \varphi(\bar b,
\bar a_{t_{\pi(0)}},\dotsc,\bar a_{t_{\pi(n-1)}})$.  Let $s_m \in I_1$
be pairwise distinct for $m < \omega$ so \nl
$\{\varphi(\bar y,\bar a_{s_n},\dotsc,\bar a_{kn+n-1}):
k < \omega\}$ is an independent contradiction.  [Saharon - details] \nl
2) Follows. \nl
3) Check directly.  \hfill$\square_{\scite{np4.4}}$
\enddemo
\bigskip

\proclaim{\stag{np4.5} Claim}  1) Assume $\bar{\bold a}^1,\bar{\bold a}^2$
are as in \scite{np4.4}.  If $\bar{\bold a}^1,\bar{\bold a}^2$ has cofinality
$> |T|$ and are mutually indiscernible and $\bar b \in {}^{\omega >} 
{\frak C}$, \ub{then} for some end-segments $J_1,J_2$ of 
Dom$(\bar{\bold a}^1),\text{Dom}(\bar{\bold a}^2)$ respectively
$\bar{\bold a}^1 \restriction J_1,\bar{\bold a}^2 \restriction J_2$ 
are mutually indiscernible over $\bar b$.
\nl
2)  Assume $\bar{\bold a}^1,\bar{\bold a}^2,\bar{\bold a}^3$ are endless
indiscernible sequences and $I = \text{ Dom}(\bar{\bold a}^\ell)$ and
$\bar a^\ell_t$ realizes {\rm Av\/}$(\{A^k_s:s <_I t \and k \in
\{1,2,3\} \cup s = t \and k < \ell\},\bar{\bold a}^\ell)$, \ub{then}:
\mr
\item "{$(a)$}"  $\langle \bar a^1_t \bar a^2_t \bar a^3_t:t \in I \rangle$ is
an indiscernible sequence;
\sn
\item "{$(b)$}"  if $\bar{\bold b}_1$ is an indiscernible set \ub{then}
$\bar{\bold b}_1,\bar{\bold b}_2$ are mutually indiscernible
\sn
\item "{$(c)$}"  if any two of $\bar{\bold a}^1,\bar{\bold a}^2,\bar{\bold a}^3$
are mutually indiscernible and $I_1,I_2,I_3$ are disjoint
unbounded subsets of $I$, \ub{then} $\bar{\bold a}^1 \restriction I_1,
\bar{\bold a}^3 \restriction I_3$ are mutually indiscernible over 
$\bar{\bold a}^2 \restriction I_2$.
\endroster
\endproclaim
\bigskip

\demo{Proof}   Similar to \scite{np4.4}(1). \nl
1) Assume not and let $I_\ell = \text{ Dom}(\bar{\bold a}^\ell)$.  We can choose
by induction on $\zeta < |T|^+$ a tuple $(n_\zeta,u^\zeta_0,u^\zeta_1,
u^\zeta_2,u^\zeta_3)$ such that:
\mr
\item "{$(*)_2(a)$}"  $|u^\zeta_\ell| = n_\zeta$
\sn
\item "{$(b)$}"  $u^\zeta_0 \cup u^\zeta_1 \subseteq I_1$
\sn
\item "{$(c)$}"  $u^\zeta_2 \cup u^\zeta_3 \subseteq I_2$
\sn
\item "{$(d)$}"  letting $\bar a^{\zeta,\ell}$ be $\langle \bar a^\zeta_t:t
\in u^\zeta_\ell \rangle$ we have $\varphi_\zeta[\bar a^{\zeta,0},
\bar a^{\zeta,2},\bar b] \and \neg \varphi_\zeta[\bar a^{\zeta,1},
\bar a^{\zeta,3},\bar b]$
\sn
\item "{$(e)$}"  $\ell \in \{0,1\} \and
\varepsilon < \zeta \and s \in u^\varepsilon_{2 \ell}
\cup u^\varepsilon_{2 \ell +1} \and t \in u^\zeta_{2 \ell} \cup 
u^\zeta_{2 \ell +1} \Rightarrow s <_{I_\ell} t$.
\ermn
Now \wilog \, $n_\zeta = n^*,\varphi_\zeta = \varphi$ for $\zeta < |T|^+$.
As $\varphi[\bar a^{\zeta,0},\bar a^{\zeta,3},\bar b]$ or $\neg \varphi
[\bar a^{\zeta,0},\bar a^{\zeta,3},\bar b]$ \wilog \, $u^\zeta_0 = 
u^\zeta_1$ or $u^\zeta_2 = u^\zeta_3$, and by the symmetry without loss of
generality the former.
Now for every $\eta \in {}^{|T|^+}2$ there is an elementary mapping
$f_\eta,f_\eta \restriction \bar{\bold a}^1$ the identity, $f_\eta$ maps
$\bar a^{\zeta,2}$ to $\bar a^{\zeta,2+\eta(\zeta)}$.  Let $g_\eta$ be an
automorphism of ${\frak C}$ extending $f^{-1}_\eta$ and let $\bar b_\eta =
g_\eta(\bar b)$.  So $\varphi[\bar a^{\zeta,0},\bar a^{\zeta,2},\bar b_\eta]$ holds
iff $\eta(\zeta)=0$ so we are done having gotten a contradiction. \nl
2) Without loss of generality $I$ is dense with no complete interval
and every interval has cardinality $> |T|$.  Now
\sn
\ub{Clause $(a)$}:

Easy.
\mn
\ub{Clause $(b)$}:

For any $s_1 <_I < \ldots <_I s_{n-1}$, by the construction we know
that: stipulating $s_0 = - \infty,s_n = + \infty$; $I_\ell = \{t \in
I:s_\ell <_I t \le_I s_{\ell +1}\}$, that the sequences $\bar{\bold
a}^1 \restriction I_0,\dotsc,\bar{\bold a}^1 \restriction I_{n-1}$ are
mutually indiscernible over $\bar a^2_{t_i} \char 94 \ldots \bar
a^2_{t_{m-1}}$. By \scite{3.4} clause (c) we are done.
\mn
\ub{Clause $(c)$}:
\mr
\item "{$(*)$}"  if $(I_1,I_2)$ are infinite disjoint intervals of $I$ then
$\bar{\bold a}^2 \restriction I_2,\bar{\bold a}^3 \restriction I_2$ are 
mutually indiscernible over $\bar{\bold a}^1 \restriction I_1 \cup
\{a^\ell_t:\ell=1,2,3 \text{ and } t \in I \backslash I_1 \backslash I_2\}$ 
\nl
[why?  by part (1), which we have already proved and (2)(a), i.e.
the indiscernibility of $\langle
\bar a^1_t \bar a^2_t \bar a^3_t:t \in I \rangle$].
\ermn
Well $(*)$ holds under any permutation of $\{1,2,3\}$ so by the way
the $a^\ell_t$'s were chosen clearly we are done.
\hfill$\square_{\scite{np4.5}}$
\enddemo
\bigskip

\proclaim{\stag{np4.6} Claim}  1) If $\bar{\bold a},\bar{\bold b}$ are
endless indiscernible, perpendicular, \ub{then} for any $\varphi(\bar
x,y,\bar c)$ for some truth values $\bold t$ we have:
\mr
\item "{$(a)$}"  for every large enough $s \in \text{ Dom}(\bar{\bold
a})$, for every large enough $t \in \text{ Dom}(\bar{\bold b})$ we
have ${\frak C} \models \varphi(\bar a_s,\bar b_t,\bar c]$
\sn
\item "{$(b)$}"  for every large enough $t \in \text{ Dom}(\bar{\bold
a})$ for every large enough $s \in \text{ Dom}(\bar{\bold a})$ we have
${\frak C} \models \varphi[\bar{\bold a}_s,\bar b_t,\bar c]$.
\endroster
\endproclaim
\bigskip

\demo{Proof}  By \scite{np4.7}(2).
\enddemo
\bigskip

\proclaim{\stag{np4.7} Claim}  1) The parallel of \scite{np4.5} holds for
several indiscernible sequences, that is, assuming $\bar{\bold a}^\zeta = 
\langle a^\zeta_t:t \in I_\zeta \rangle$ is an endless indiscernible sequence
for $\zeta < \zeta^*$
\mr
\item "{$(A)$}"  If the intervals $[\text{cf}(I_\zeta),|I_\zeta|]$ are 
pairwise disjoint, cf$(I_\zeta) > |T| + \zeta^*$, \ub{then} for some end
segment $J_\zeta$ of $I_\zeta$ for $\zeta < \zeta^*$, we have $\langle
\bar{\bold a}^\zeta \restriction J_\zeta:\zeta < \zeta^* \rangle$ is mutually
indiscernible, which means: each $\bar{\bold a}^\zeta \restriction
\bold J_\zeta$ is indiscernible over $\cup\{\bar{\bold a}^\varepsilon
\restriction J_\varepsilon:\varepsilon < \zeta^* \and \varepsilon \ne \zeta\}$
(in fact we can get indiscernibility over $\cup\{\bar{\bold a}^\varepsilon:
\varepsilon < \zeta^* \and \varepsilon \ne \zeta\})$
\sn
\item "{$(B)$}"  Assume $\langle \bar{\bold a}^\zeta:\zeta < \zeta^* \rangle$
are mutually indiscernible, $\bar b \in {}^{\omega >}{\frak C}$ and $I_\zeta =
\text{ Dom}(\bar{\bold a}^\zeta)$ and $I_\zeta = \text{ cf}(\text{Dom}
(\bar{\bold a}^\zeta)) > |T| + \zeta^*$.
\ub{Then} there are end segments $J_\zeta$ of $I_\zeta$ for $\zeta < \zeta^*$
such that $\langle \bar{\bold a}^\zeta \restriction I_\zeta:\zeta < \zeta^*
\rangle$ is mutually indiscernible over $\bar b$.
\sn
\item "{$(C)$}"  If $J$ is an infinite linear order disjoint to $\cup
\{I_\zeta:\zeta < \zeta^*\}$ and $a^\zeta_t$ realizes Av$(\{\bar
a^\varepsilon_s:\varepsilon < \zeta^*$ and $s \in I_\varepsilon$
\ub{or} $s \in J \and t <_J s$ \ub{or} $s = t \and \varepsilon <
\zeta\},\bar{\bold a}^\zeta)$ \ub{then} $\{\langle \bar
a^{\varepsilon^\zeta}_s:s \in J \rangle:\zeta < \zeta^* \rangle$ are
mutually indiscernible over $\cup\{\bar a^\varepsilon_s:\varepsilon <
\zeta^*,s \in I_\varepsilon\}$.
\ermn
2) We weaken in the conclusion the mutually indiscernible by mutually
$\Delta$-indiscernible, \ub{then} we can weaken cf$(I_\zeta) > |T| +
|\zeta^*|$ to cf$(I_\zeta) > |\zeta^*|$.
\endproclaim
\bigskip

\demo{Proof} Easy.
\enddemo
\bigskip

\proclaim{\stag{np4.8} Claim}  1) If $\bar{\bold a} = \langle \bar a_t:t \in I
\rangle$ is an indiscernible sequence, $\bar b \in {}^{\omega >}{\frak C}$ 
\ub{then} we can divide $I$ to $\le 2^{|T|}$ convex subsets $\langle I_\zeta:
\zeta < \zeta^* \rangle$ such that $\langle \bar{\bold a} \restriction 
I_\zeta:\zeta < \zeta^*,I_\zeta \text{ infinite}\rangle$ is mutually 
indiscernible over $\bar b$. \nl
2) Similarly in \scite{np4.7}.
\endproclaim
\bigskip

\demo{Proof}  Easy.
\enddemo
\bigskip

\proclaim{\stag{np4.8a} Claim}  Assume that
\mr
\item "{$(a)$}"  $\bar{\bold b},\bar{\bold a}^\zeta$ are an endless
indiscernible sequence for $\zeta < \zeta^*$
\sn
\item "{$(b)$}"  $\bar{\bold a}^\zeta,\bar{\bold a}^\varepsilon$ are
perpendicular for $\zeta \ne \varepsilon$
\sn
\item "{$(c)$}"  $\bar{\bold b},\bar{\bold a}^\zeta$ are not perpendicular.
\ermn
\ub{Then} $\zeta^* < |T|^+$.
\endproclaim
\bigskip

\demo{Proof}  Assume toward contradiction that $\zeta^* \ge |T|^+$.  We let
$A = \bar{\bold b} \cup \dbcu_\zeta \bar{\bold a}^\zeta$ and choose
$\bar a^{\zeta,*}_n,\bar b^*_n$ for $n < \omega$ such that:
\mr
\item "{$(a)$}"  $\bar a^{\zeta,*}_n$ realized Av$(A \cup\{b^*_m:m < n\}
\cup \{a^{\varepsilon,*}_m:m < n \and \varepsilon < \zeta^*$ or $m = n \and
\varepsilon < \zeta\},\bar{\bold a}^\zeta)$
\sn
\item "{$(b)$}"  $\bar b^*_n$ realizes Av$(A \cup \{b^*_m:m < n\} \cup
\{a^{\varepsilon,*}_m:m \le n,\varepsilon < \zeta^*\},\bar{\bold b})$.
\ermn
For each $\zeta$, as $\bar{\bold b},\bar{\bold a}^\zeta$ are not 
perpendicular, we can find $n_\zeta < \omega,u^\ell_\zeta \in [\omega]
^{n_\zeta}$ for $\ell = 0,1,2$ such that $\langle \bar b^*_n:n \in
u^0_\zeta \rangle \char 94 \langle \bar a^{\zeta,*}_n:n \in u^1_\zeta \rangle$
and $\langle \bar b^*_n:n \in u^0_\zeta \rangle \char 94 \langle
\bar a^{\zeta,*}_n:n \in u^2_\zeta \rangle$ does not realize the same type;
say one satisfies $\varphi_\zeta(\bar x,\bar y)$ the second not. As we can
replace $\langle \bar{\bold a}^\zeta:\zeta < |T|^+ \rangle$ by any
subsequence of length $|T|^+$, without loss of generality $n_\zeta = n_*,
u^\ell_\zeta = u_\ell,\varphi_\zeta = \varphi$.  Now for every ${\Cal U}
\subseteq |T|^+$ let $f_{\Cal U}$ be the elementary mapping with domain
$\cup\{a^{\zeta,*}_n:n \in u_1,\zeta < |T|^+\}$, mapping $a^{\zeta,*}_{n_1}$
to $a^{\zeta,*}_{n_2}$ iff $\zeta \in {\Cal U},n_1=n_2$ or $\zeta \in
|T|^+ \backslash {\Cal U},n_1 \in u_1,n_2 \in u_2,|n_2 \cap u_1| =
|n_2 \cap u_2|$.  Let $g_{{\Cal U}}$ be an automorphism of ${\frak C}$ extending
$f^{-1}_{\Cal U}$.  We have gotten the independence property for $\varphi
(\bar x,\bar y)$ as $g_{\Cal U}(\langle \bar b^*_n:n \in u^0_0 \rangle)$
realizes $\{\varphi(\langle \bar x_n:n \in u_0 \rangle,\langle 
\bar b^{\zeta,*}_n:n \in u_1 \rangle)^{\text{if}(\zeta \in {\Cal U})}:\zeta
< |T|^+\}$, contradiction.  \hfill$\square_{\scite{np4.8a}}$
\enddemo
\bn
\centerline {$* \qquad * \qquad *$}
\bn
Recall (\cite[Ch.III,\S7]{Sh:c})
\definition{\stag{np4.9} Definition}  1) $p \in \bold F^{sp}_\kappa(B)$ if
for some set $A$ we have
$p \in S^{< \omega}(A),B \subseteq A,|B| < \kappa$ and $p$ does not split
over $B$. \nl
2) ${\Cal A} = (A,\langle \bar b_i,B_i:i < i^* \rangle)$ is an
$\bold F^{sp}_\kappa$-construction (or $\langle b_i,B_i:i < i^* \rangle$
is an $\bold F^{sp}_\kappa$-construction over $A$) \ub{if} tp$(\bar b_i,A \cup
\{b_j:j < i\}) \in F^{sp}_\kappa(B_i)$, so $B_i \subseteq A^{\Cal A}_i =:
A \cup \{b_j:j < i\})$. \nl
3) Omitting $B_i$ means for some $B_i$; let $i^* = \ell g({\Cal A})$.
\enddefinition
\bigskip

\remark{Remark}  We may use $\bar b_i$'s of any length $< \kappa$.
\endremark
\bigskip

\proclaim{\stag{np4.10} Claim}  1) If $B \subseteq A,p$ is an $m$-type over
$B$, \ub{then} there is $q \in S^m(A)$ extending $p$ and $B_1 \subseteq A,
|B| \le |T|$ such that $q$ does not split over $B \cup B_1$. \nl
2) For any $A$ and $\kappa = \text{ cf}(\kappa) > |T|$ there is a model $M$
and $\bold F^{sp}_\kappa$-construction ${\Cal A} = (A,\langle \bar b_i,B_i:i < i^*
\rangle)$ such that:
\mr
\item "{$(a)$}"  $M = A^{\Cal A}_{i^*},\|M\| = |A|^{< \kappa} +
\dsize \sum_{\theta < \kappa} 2^{2^\theta}$
\sn
\item "{$(b)$}"  $M$ is $\kappa$-saturated, moreover if $B \subseteq M,
|B| < \kappa,p \in S^m(M)$ does not split over $B$ then for unboundedly
many $i < i^*,\bar b_\ell$ realizes $p \restriction A^{\Cal A}_i$
\sn
\item "{$(c)$}"  cf$(i^*) \ge \kappa$
\ermn
3) If ${\Cal A}$ is an $\bold F^{sp}_\kappa$-construction, $\kappa =
\text{ cf}(\kappa),\bar b \subseteq A^{\Cal A}_{\ell g({\Cal A})}$ has
length $< \kappa$, \ub{then} tp$(\bar b,A)$ does not split over some
$B \subseteq A,|B| < \kappa$.
\endproclaim
\newpage

\head {\S5 Indiscernible sequence perpendicular to cuts} \endhead  \resetall \sectno=5
\bigskip

Our aim is to show that for a set of $\{\bar{\bold b}_\zeta:\zeta <
\zeta^*\}$ of pairwise perpendicular endless indiscernible sets, we
can find a model $M \supseteq \cup \{\bar{\bold b}_\zeta:\zeta <
\zeta^*\}$ with $\langle \text{dual-cf}(\bar{\bold b}_\zeta):\zeta <
\zeta^* \rangle$ essentially as we like, and other $\bar{\bold b}'$ in
$M$ has such dual cofinality iff this essentially follows.  Toward
this we define and investigate when an endless indiscernible sequence
$\bar c$ is perpendicular to a (Dedekind) cut $(I_1,I_2)$ is an
indiscernible sequence $\bar{\bold a}$.
\bigskip

\definition{\stag{np5.1} Definition}  1) We say $(I_1,I_2)$ is a Dedekind cut
of the linear order $I$, if $I$ is the disjoint union of $I_1,I_2$ and
$s \in I_1 \and t \in I_2 \Rightarrow s <_I t$ and we write
$I,I = I_1 + I_2$, and its cofinality is $(\text{cf}(I_1),\text{cf}
(I^*_2))$.  If $I$ is a convex subset of $J$ and
$I_1 \ne \emptyset \ne I_2$ we may abuse our notation saying
``$(I_1,I_2)$ is a Dedekind cut of $J$".
We say
$(I_1,I_2)$ is a Dedekind cut of $\bar{\bold a}$ if it is a Dedekind cut of
Dom$(\bar{\bold a})$.  If not say otherwise, $I_1 \ne \emptyset \ne I_2$,
and the cut is nontrivial if both its cofinalities are infinite. \nl
2)  $(J_1,J_2) \le (I_1,I_2)$ if $J_1$ is an end segment of $I_1$ and $J_2$
is an initial segment of $I_2$. \nl
3) We say the set $A$ respects the Dedekind cut $(I_1,I_2)$ of $\bar{\bold a}$
\ub{if} $(I_1,I_2)$ is a Dedekind cut of $\bar{\bold a}$ and for every
$\bar b \in {}^{\omega >}A$ for some $(J_1,J_2) \le (I_1,I_2)$ the sequence
$\bar{\bold a} \restriction (J_1 + J_2)$ is indiscernible over $\bar b$. \nl
4) For endless indiscernible sequences 
$\bar{\bold a},\bar{\bold b}$ and Dedekind cut $(I_1,I_2)$ of $\bar{\bold a}$ we say that
$\bar{\bold b}$ is perpendicular to the cut when: 
if $\bar{\bold b}'$ is an indiscernible
sequence over $\bar{\bold b} \cup \bar{\bold a}$ based on $\bar{\bold
b}$ (see below) \ub{then} $\bar{\bold b}' \cup \bar{\bold a}$ respects the cut 
$(I_1,I_2)$ of $\bar{\bold a}$. \nl
5) For endless indiscernible sequences $\bar{\bold a}$ and $A \supseteq
\bar{\bold a}$ we say an endless indiscernible sequence $\bar{\bold b} = \langle
\bar{\bold b}_i:t \in I \rangle$ over $A$ is based on $\bar{\bold a}$ 
if each $\bar b_t$ realizes Av$(A \cup \{\bar b_s:s <_I t\},\bar{\bold a})$.
\enddefinition
\bigskip

\proclaim{\stag{np5.2} Claim}  1) If $\langle A_i:i < \delta \rangle$ is
increasing, $\bar{\bold a} \subseteq A_0$ is an endless indiscernible
sequence, $\bar a'_i \subseteq A_{i+1}$ realizes {\rm Av\/}$(A_i,\bar{\bold a}),
\bar{\bold a}' = \langle \bar a'_i:i < \delta \rangle,\bar{\bold a}''$ is
the inverse of $\bar{\bold a}'$ \ub{then}
\mr
\item "{$(a)$}"  $\bar{\bold a} \char 94 \bar{\bold a}''$ is indiscernible
\sn
\item "{$(b)$}"  the set $\dbcu_{i < \delta} A_i$ respects the cut
$(\text{Dom}(\bar{\bold a}),\text{Dom}(\bar{\bold a}''))$ of
$\bar{\bold a} \char 94 \bar{\bold a}''$.
\ermn
2) If $\bar{\bold a}$ is a non stable indiscernible sequence, 
$\bar{\bold a} \subseteq A$,
the set $A$ respects the endless cut $(I_1,I_2)$ of $\bar{\bold a}$ and the 
cofinalities of
the cut are $> |T|$ \ub{then} {\rm dual-cf\/}$(\bar{\bold a} \restriction I_1,M) =
\text{ cf}(I^*_2)$. \nl
3) If $\bar{\bold a}$ is an indiscernible sequence with Dedekind cut
$(I_1,I_2)$ of cofinaltiy $(\kappa_1,\kappa_2),\aleph_0 \le
\kappa_1,\kappa_2$ and $\bar{\bold c}$ an endless indiscernible
sequence respecting this cut then: for some for every formula
$\varphi(\bar x,\bar y,\bar z)$ and sequence $\bar b$ for some truth
value $\bold t$ we have:
\mr
\item "{$(*)(i)$}"  for every large enough $s \in \text{
Dom}(\bar{\bold c})$, for some $(J_1,J_2) \le (I_1,I_2)$ for every $t
\in J_1 \cup J_2$ we have ${\frak C} \models \varphi[\bar a_2,\bar
b,\bar c_1]^{\bold t}$,
\sn
\item "{${{}}(ii)$}"  for some $(J_1,J_2) \le (I_1,I_2)$ for every $t
\in J_1 \cup J_2$ we have ${\frak C} \models \varphi[\bar a_2,\bar
b,\bar c_1]^{\bold t}$, ???? and end segment $J$.
\ermn
4) If in part (3), $|A| + |A| < \kappa_1,\kappa_2$ then for some
$(J_1,J_2) \le (I_1,I_2)$ we have: of Dom$(\bar{\bold c}),\bar{\bold
a} \restriction (J_1 + J_2),\bar{\bold c} \restriction J$ are mutually
$\Delta$-indiscernible.
\endproclaim
\bigskip

\demo{Proof}  1), 2) Straightforward.
\nl
3) Let $\delta = |T|^+,\bar c_\gamma$ realizes Av$(\bar{\bold a} \cup
\bar{\bold c} \cup \{\bar c_\beta:\beta < \gamma\},\bar{\bold c})$,
for $\gamma < \delta$ so by the definition of ``respect the Dedekind
cut", there is $(J_1,J_2) \le (I_1,I_2)$ such that $\bar{\bold a}
\restriction (J_1 \cup J_2),\langle \bar c_\gamma:\gamma < \delta
\rangle$ are mutually indiscernible.  Let $(I_1,I_2)$ have cofinality
$\kappa_1,\kappa_2)$ and for our purpose \wilog \, $\kappa_1,\kappa_2
> |T|$.  Now $\bar{\bold a} \restriction J_1$, the inverse of
$\bar{\bold a} \restriction J_2,\langle \bar c_\gamma:\gamma < \delta
\rangle$ are mutually indiscernible, hence by \scite{np4.7}, clause
(B) \wilog \, they are mutually indiscernible over $\bar b$
(i.e. omitting an initial segment of each and renaming.  So we hafe
truth values $\bold t(1),\bold t(2)$ such that $t \in J_\ell \and
\gamma < \delta \Rightarrow {\frak C} \models \varphi[\bar a_t,\bar
b,\bar c_\gamma]^{{\bold t}(\ell)}$.  If $\bold t(1) \ne \bold t(2)$
we get contradiction to ``$T$ has NIP" so \wilog \, $\bold t(1) =
t(2)$ and as we can replace $\varphi$ by $\neg \varphi$ \wilog \,
$\bold t(1) = \bold t(2) =$ truth.  So by the choice of $\langle \bar
c_\gamma:\gamma < \delta \rangle$, for every $t \in J_1 \cup J_2$, for
every large enough $s \in \text{ Dom}(\bar{\bold c})$ we have ${\frak
C} \models \varphi[\bar a_t,\bar b,\bar c_s]$.

Clearly $\bar{\bold c},\bar{\bold a} \restriction J_1$ is
perpendicular (by \scite{np4.6}).
\enddemo
\bigskip

\proclaim{\stag{np5.3} Claim}  Assume
\mr
\item "{$(a)$}"  $I = I_1 + I_2$ and the Dedekind cut $(I_1,I_2)$ has
cofinality $(\kappa_1,\kappa_2)$
\sn
\item "{$(b)$}"  $\bar{\bold a} = \langle \bar a_t:t \in I \rangle$ is
an indiscernible sequence
\sn
\item "{$(c)$}"  $\bar{\bold a} \subseteq A$
\sn
\item "{$(d)$}"  the set $A$ respects the cut $(I_1,I_2)$ of $\bar{\bold a}$.
\ermn
1) If tp$(\bar b,A) \in \bold F^{sp}_\kappa$ and $\kappa \le \kappa_1,
\kappa_2$, \ub{then} the set $A \cup \bar b$ respects the cut $(I_1,I_2)$
of $\bar{\bold a}$. \nl
Assume in addition
\mr
\item "{$(e)$}"  $|T| < \kappa_1,\kappa_2$
{\roster
\itemitem{ $(A)$ }  If 
$\bar{\bold c} \subseteq A$ is an endless indiscernible sequence and
$\bar{\bold c}$ is perpendicular to the cut $(I_1,I_2)$ of $\bar{\bold a}$
and $\bar c$ realizes {\rm Av\/}$(A,\bar{\bold c})$, \ub{then} $A \cup \bar c$
respects the cut $(I_1,I_2)$ of $\bar{\bold a}$. 
\sn
\itemitem{ $(B)$ }  If $A_i \, (i < \delta)$ is increasing each $A_i$ respects the cut
$(I_1,I_2)$ of $\bar{\bold a}$ then also $\dbcu_{i < \delta} A_i$
does. 
\endroster}
\ermn
2) If $A^+ = A \cup \{a_i:i < i^* \}$ and for each $i$, tp$(a_i,A \cup
\{a_j:i < i^*\})$ belongs to $\bold F^{sp}_{\text{min}\{\kappa_1,\kappa_2\}}$
\ub{or} is {\rm Av\/}$(A \cup \{a_j:j < i\},\bar{\bold b})$ where $\bar{\bold b} 
\subseteq A \cup \{a_j:j < i\}$ is an endless indiscernible sequence 
perpendicular to the cut
$(I_1,I_2)$ of $\bar{\bold a}$, \ub{then} $A^+$ respects the cut $(I_1,I_2)$
of $\bar{\bold a}$.
\endproclaim
\bigskip

\demo{Proof}  1) Check. \nl
2) Suppose that this fails, so there is $b \in {}^{\omega >}(A)$ such
that $\bar b \char 94 \bar c$ witness it.  Now by assume (e) for some
$(J_1,J_2) \le (I_1,I_2)$ we have $\bar a \restriction J_1,\bar{\bold
a} \restriction J_2$ are mutually indiscernible over $\bar b \char 94
\bar c$.  As $\bar b \char 94 \bar c$ witness failure for some
$\varphi(\bar x,\bar y,\bar z)$ we have
\mr
\item "{$(*)_1(\alpha)$}"  $t \in J_1 \Rightarrow {\frak C} \models
\varphi[\bar a_2,\bar b,\bar c]$ and
\sn
\item "{${{}}(\beta)$}"  $t \in J_2 \Rightarrow {\frak C} \models \neg
\varphi[\bar a_t,\bar b,\bar c]$.
\ermn
By clause $(\alpha)$, for every $t \in J_1$ for every large enough $s
\in \text{ Dom}(\bar{\bold c})$ we have ${\frak C} \models
\varphi[\bar a_t,\bar b,\bar c_s]$, however, $\bar c$ is perpendicular
to the cut $(I_1,I_2)$ of $\bar{\bold a}$ hence for every large enough
$s \in \text{ Dom}(\bar{\bold c})$ for every large enough $t_1 \in
J_1$ and small enough $t_2 \in J_2$ we have ${\frak C} \models
\varphi[\bar a_{t_1},\bar b,\bar c_s] \and \varphi[\bar a_{t_2},\bar
b,\bar c_s]$.  Again as $\bar c$ is perpendicualr to the cut
$(I_1,I_2)$ of $\bar{\bold a}$ we get: for every small enough $t_2 \in
J_2$ for every large enough $s \in \text{ Dom}(\bar{\bold c})$ we have
${\frak C} \models \varphi[\bar a_{t_2},\bar b,\bar c_1]$,
contradicting clause $(\beta)$ of $(*)$.  \nl
3) Check the definition. \nl
4) Prove by induction on $i$ using (1), (2), (3).
\hfill$\square_{\scite{np5.3}}$
\enddemo
\bigskip

\proclaim{\stag{np5.4} Claim} Assume
\mr
\item "{$(a)$}"  $(I_1,I_2)$ is a cut of the indiscernible sequence
$\bar{\bold a}$ with both cofinalities infinite
\sn
\item "{$(b)$}"  $\bar{\bold b}$ is an endless indiscernible sequence
\sn
\item "{$(c)$}"  $\bar{\bold a} \restriction I_1,\bar{\bold b}$ are
perpendicular
\sn
\item "{$(d)$}"  for $t \in I_2,\bar a^1_t$ realizes {\rm
Av\/}$(\{\bar a^1_s:s \in I_1
\vee t <_I s \in I_2\} \cup \bar{\bold b},\bar{\bold a}^1 \restriction
I_1)$.
\ermn
\ub{Then} $\bar{\bold b}$ is perpendicular to the cut $(I_1,I_2)$ of
$\bar{\bold a}$.
\endproclaim
\bigskip

\demo{Proof}  First assume that the cofinalities of $I_1,I^*_2,\bar{\bold b}$
are $> |T|$.  
By the demands (c) + (d), $\bar{\bold a} \cup \bar{\bold b}$ respect
the cut $I_1,I_2)$ of $\bar{\bold a}$.
So assume $A \supseteq \bar{\bold a} \cup \bar{\bold b}$ 
respect the cut $(I_1,I_2)$ of $\bar{\bold a}$ and $\bar b$ realizes
Av$(A,\bar{\bold b})$.  So let $\bar c \in {}^{\omega >}A$; as $A$ respects
the cut $(I_1,I_2)$ of $\bar{\bold a}$, there is $(J_1,J_2) \le (I_1,I_2)$
such that $\bar{\bold a} \restriction (J_1 + J_2)$ is indiscernible over
$\bar c$ and no problem.

Generally deal with $\Delta$-types for finite $\Delta$'s.
\hfill$\square_{\scite{np5.4}}$
\enddemo
\bigskip

\proclaim{\stag{np5.5} Claim}  Assume
\mr
\item "{$(a)$}"  $\lambda = \lambda^{< \kappa_2}$
\sn
\item "{$(b)$}"  $\kappa_1 = \text{ cf}(\kappa_1) \le \kappa_2 \le
\theta_2 = \text{ cf}(\theta_2),\kappa_1 \le \theta_1 = \text{ cf}(\theta_1)
\le \lambda$
\sn
\item "{$(c)$}"  $|A| \le \lambda$
\sn
\item "{$(d)$}"  $\bar{\bold a}^\zeta \subseteq A$ is endless, nonstable
indiscernible for $\zeta < \zeta^*$ and $\zeta^* \le \lambda$
\sn
\item "{$(e)$}"  the $\bar{\bold a}^\zeta$ for $\zeta < \zeta^*$ are
pairwise perpendicular.
\ermn
\ub{Then} we can find a model $M$ such that
\mr
\item "{$(\alpha)$}"  $A \subseteq M$
\sn
\item "{$(\beta)$}"  {\rm dual-cf\/}$(\bar{\bold a}^\zeta,M) = \theta_1$ for
every $\zeta < \zeta^*$
\sn
\item "{$(\gamma)$}"  if $\bar a \subseteq M$ is a nonstable endless
indiscernible sequence of cardinality (hence cofinality) $< \kappa_2$
perpendicular to every $\bar{\bold a}^\zeta$ \ub{then} 
{\rm dual-cf\/}$(\bar{\bold a},M) = \theta_2$
\sn
\item "{$(\delta)$}"  $M$ is $\kappa_1$-saturated.
\endroster
\endproclaim
\bigskip

\demo{Proof}  We first deal with a restricted case, then derive from it the
general case.
\enddemo
\bn
\ub{Case 1}:  $|A| < \lambda = \text{ cf}(\lambda) = \theta_1,\kappa_1 =
\kappa_2 = \text{ cf}(\kappa_2) \le \theta_2$ and $(\forall \alpha < \lambda)
(|\alpha|^{\theta_2} < \lambda)$ and in clause $(\gamma)$ we demand just
dual-cf$(\bar{\bold a},M) \in [\theta_2,\lambda)$.

We can find $a_i$ for $i < \lambda$ such that letting $A_i = A \cup
\{\bar a_j:j < i\}$ we have
\mr
\item "{$(i)$}"  for each $i < \lambda$ we have tp$(\bar a_i,A_i) \in
\bold F^{sp}_{\theta_2}$ or tp$(a_i,A_i) =$ {\rm Av\/}$(A_i,\bar{\bold a}^\zeta)$ for
some $\zeta < \zeta^*$
\sn
\item "{$(ii)$}"  if $p \in S^{< \omega}(A_\lambda),p \in \bold F^{sp}
_{\theta_2}$ \ub{then} for $\lambda$ ordinals $j  < \lambda,p$ is realized
by $b_j$.
\ermn
This is straightforward and clauses $(\alpha),(\beta),(\delta)$ obviously
hold.  As for clause $(\gamma)^-$, let $\bar{\bold a} = \langle a_t:t \in I
\rangle \subseteq M$ endless indiscernible $|I| < \theta_2,\bar{\bold a}$
perpendicular to every $\bar{\bold a}^\zeta$ and assume toward contradiction
that dual-cf$(\bar{\bold a},M) \notin [\theta_2,\lambda)$ but trivially it
is $\le \|M\| = \lambda$ and by saturation $\ge \theta_2$, so necessarily
it is $\lambda$. \nl
So for every $\alpha < \lambda$ some $\bar c_\alpha \in M$ realizes
Av$(A_\alpha,\bar{\bold a})$ and let $\beta_\alpha = \text{ Min}\{\beta <
\lambda:\bar c_\alpha \subseteq A_\beta\}$.  So $\beta_\alpha \in (\alpha,
\lambda)$ and let $E = \{\delta < \lambda:\delta$ limit $\bar{\bold a}
\subseteq A_\delta$ and $(\forall \alpha < \delta),\beta_\alpha < \delta\}$.
For $\delta \in \text{ acc}(E)$ let $\bar a^{\delta,0}$ be the sequence
$\langle \bar c_\alpha:\alpha \in E \cap \delta \rangle$ and 
$\bar{\bold a}^\delta$ be its inverse.  Now not only is $\bar{\bold a}
\char 94 \bar{\bold a}^\delta$ is indiscernible but $A_\delta$ respects the
cut $(\text{Dom}(\bar{\bold a}),\text{Dom}(\bar{\bold a}^\delta))$ of
$\bar{\bold a}$ to \scite{np5.2}(2).  Choose $\delta(*) \in \text{ acc}(E)$
be of cofinality $\theta_2$.  Now by \scite{np5.3}(4) $\bar{\bold a}
\char 94 \bar{\bold a}^{\delta(*)}$ is respected also by $A_\lambda = M$.
By \scite{np5.2}(2) this gives dual-cf$(\bar{\bold a},M) = \text{ cf}
(\delta(*)) = \theta_2$, contradiction.
\bn
\ub{Case 2}:  As above getting the full $(\gamma)$.  
By absoluteness arguments without loss of generality
we can find $\theta_*,\lambda_*$ such
that $\theta_* = \theta^{< \theta_*}_*,2^{\theta_*} = \lambda_*$ and they are
$> \lambda$.  Now use case A for $\theta_*,\lambda_*$ getting $M_0$. We 
can find $M$ such that
\mr
\item "{$(*)(a)$}" $M_1 \prec M_0$ include $A$ and has cardinality
$\theta_*,\kappa$-saturated
\sn
\item "{$(b)$}"  for $\zeta < \zeta^*$, dual-cf$(\bar{\bold a}^\zeta,M_1)
= \theta_1$
\sn
\item "{$(c)$}"  any endless nonstable indiscernible sequence $\bar{\bold a}
\subseteq M_1,|\text{Dom}(\bar{\bold a})| < \theta_2$ perpendicular to every
$\bar{\bold a}^\zeta$, dim-cf$(\bar{\bold a}) = \theta_2$ \nl
[why it exists?  we choose by induction  on $\alpha < \theta_1,M_{1,\alpha}$
satisfying clause $(*)(a)$, increasing continuous with $\alpha$ such that
if $\bar{\bold a} \subseteq M_{1,\alpha}$ is as in (2), then a witness to
dim-cf$(\bar{\bold a},M) = \theta_2$ is included in $M_{1,\alpha+1}$ and
Av$(M_{1,\alpha},\bar{\bold a}^\zeta)$ is realized in $M_{1,\alpha+1}$.
Now $\dbcu_{\alpha < \theta_1} M_{1,\alpha}$ is as required.
\ermn
Now similarly we can find $M_2 \prec M_1$ as required this time by a sequence
$\theta_2$ approximations.  \hfill$\square_{\scite{np5.5}}$
\newpage

\head {\S6 Concluding Remark} \endhead  \resetall \sectno=6
\bn
\ub{\stag{np6.1} Discussion}:  A major lack of this 
work is the absence of test questions. \nl
A candidate is (see \cite[\S2]{Sh:702}).
\mn
\ub{\stag{np6.2} Question}:  If $A \subseteq {\frak C}_T,\kappa = |A| + |T|$ (or
$\kappa = \beth_7(|A|+|T|)$ and $\lambda = \beth(2^\kappa)^+$ (or larger,
but \ub{no} large cardinals) and $a_i \in {\frak C}_T$ for $i < \lambda$
then for some $w \in [\lambda]^{\kappa^+}$, the sequence $\langle a_i:i \in
\omega \rangle$ is an indiscernible sequence over $A$ (in ${\frak C}_T$).
\mn
Through this property does not characterize NIP, it is quite natural in
this context.  See also next.
\mn
Another direction is generalizing DOP, which in spite of its name is
a non first order independence property.
\bigskip

\definition{\stag{np6.3} Definition}  $T$ has the dual-cf-$\kappa$-dimensional
independence if: $\bar \kappa = (\kappa_0,\kappa_1,\kappa_2)$, \nl
$\kappa_1 \ne
\kappa_2 > \kappa_0 < \kappa_1,\kappa_0 < \kappa_2$ and for every $\lambda$
and $R \subseteq \lambda \times \lambda$ symmetric we can find $M_R,
\bar{\bold b}_\alpha,\bar{\bold c}_\alpha \in {}^{\kappa_0}(M_R)$ and
$\bold I_{\alpha,\beta} = \langle \bar a_{\alpha,\beta,i}:i < \kappa_0
\rangle \subseteq M_R$ for $(\alpha,\beta) \in R,\alpha < \beta$ such that:
\mr
\item "{$(a)$}"  the type of $\bar{\bold b}_\alpha \char 94
\bar{\bold c}_\beta \char 94 \bold I_{\alpha,\beta}$ is the same for all
$\alpha,\beta$
\sn
\item "{$(b)$}"  dual-cf$(\bold I_{\alpha,\beta},M_R) = \kappa_1$
\sn
\item "{$(c)$}"  if $\alpha < \beta \neg \alpha R_\beta,\bold I'_{\alpha,\beta} = \langle
a'_{\alpha,\beta}:i < \kappa_0 \rangle \subseteq M_R$ 
is such that for every $(\alpha_1,\beta_1) \in R$ there is an
automorphism $h$ of ${\frak C}$ taking $\bar{\bold b}_{\alpha_1}$ to
$\bar{\bold b}_\alpha,\bar{\bold c}_{\beta_1}$ to $\bar{\bold
c}_\beta$ and $\bar a_{\alpha_1,\beta_1}$ to $a_{\alpha,\beta,i}$
\ub{then} dual-cf$(\bold I'_{\alpha,\beta},M) = \kappa_2$
\sn
\item "{$(d)$}"  $M_R$ is $\kappa^+_0$-saturated.
\ermn
Note that (d), (d) follows from
\mr
\item "{$(c)^+$}"  if $\bold I'_{\alpha,\beta} = \langle \bar a'_{\alpha,
\beta,i}:i < \kappa_0 \rangle$ realizes the relevant type
$(\alpha,\beta) \notin R,\alpha < \beta < \lambda$ and $\alpha_1 <
\beta_1 < \lambda,(\alpha_1,\beta_1) \in R$ then $\bold I'_{\alpha,\beta},
\bold I_{\alpha_1,\beta_1}$ are perpendicular.
\ermn
As in \S5 we can show that many variants are equivalent (using $+ \infty,
-\infty$ to absorb).  We can similarly discuss deepness.
\enddefinition 
\bn
\ub{\stag{np6.4} Discussion}:  1) It is know that e.g. the $p$-adics are NIP (but
unstable).  Does this work tell us anything on them?  Well, the construction
in \S5 gives somewhat more than what unstability gives: complicated models
with more specific freedom.  Note that instead dual-cf$(\bold I,M)$ we can
use more complicated invariants (see \cite[Ch.III,\S3]{Sh:e}) or earlier
works).

We can, of course, (for the $p$-adic) characterize directly when 
indiscernible sequences are perpendicular.

2) We may like to define super-NIP (and $\kappa_{\text{nip}}(T)$) (parallel
of superstable, i.e. $\kappa(T) = \aleph_0$ or super simple
$\kappa_{cdt}(T) = \aleph_0$).  This is not clear to me.  We may try the
definition ``$w(\bold I) < \aleph_0$" for every endless indiscernible
sequence where
\definition{\stag{np6.5} Definition}  For an endless indiscernible sequence $\bold I$
let $w(\bold I) = \sup\{\alpha:$ there is a sequene of length $\alpha$ of
pairwise endless indiscernible sequences each non perpendicular to
$\bold I\}$.  But $w(\bold I)$ is not exactly like dimension in
the sense of algebraic manifolds.
\mn
\ub{Question}:  Assume 
$\bold I_\ell = \langle a^\ell_t:t \in I_\ell \rangle$ for $\ell =1,2$
are endless indiscernible nonperpendicular sequences
\mr
\item "{$(a)$}"  find a definable equivalence relation $E$ such that 
$\langle a^2_t/E:t \in I_r \rangle$ is nontrivial and $a^2_t \in
\text{ acl}(\bold I_1 \cup\{a^2_s:s <_{I_2} t\})$ for any large enough $t$
\sn
\item "{$(b)$}"  if $(\bold I_1,\bold I_2)$ is $(1,< \omega)$-mutual
indiscernible, can we define a derived group?  More generally, it seems
persuasive that groups appear naturally, particularly ordered groups
\sn
\item "{$(c)$}"  does the fact that putting of elements together, make
strong splitting to dividing helps?
\sn
\item "{$(d)$}"  can the canonical bases of \S1 help?  Do they help for
simple theories
\sn
\item "{$(e)$}"  what can we say on ``$\bar{\bold a}$ orthogonal to a set
model $A$?"
\ermn
Can we say more on  ``stable" aspects?  (see \S1).

Cherlin wonders on the place of parallel algebraic geometric dimension and
place of 0-minimal theory.  In my perception probably if we succeed in 2),
we may have a minimality notion which may then be characterized as some
cases, but maybe it does not fit.
\bn
\ub{\stag{np6.6} Question}:  Given two non perpendicular types which are weakly
perpendicular can we find naturally defined groups?
\enddefinition
\bn
\centerline {$* \qquad * \qquad *$}
\bn
\proclaim{\stag{np8.1} Claim}  Assume
\mr
\item "{$(\alpha)$}"  $\bar{\bold b}^0 = \langle \bar b_t:t \in I_0
\rangle$ is an infinite indiscernible sequence over $A$
\sn
\item "{$(\beta)$}"  $B \subseteq {\frak C}$.
\ermn
\ub{Then} we can find $I_1$ and $\bar b_t$ for $t \in I_1 \backslash
I_0$ such that:
\mr
\item "{$(a)$}"  $I_0 \subseteq I_1,|I_1 \backslash I_0| \le |B| +
|T|$
\sn
\item "{$(b)$}"  $\bar{\bold b}' = \langle \bar b_t:t \in I_1 \rangle$
is an indiscernible sequence over $A$
\sn
\item "{$(c)$}"  if $I_1 \subseteq I_2$ and $\bar b_t$ for $t \in I_2
\backslash I_1$ are such that $\bar{\bold b}^2 = \langle \bar b_t:t \in I_1 \rangle$
is an indiscernble sequence over $A$.
\endroster
\endproclaim
\bigskip

\demo{Proof}  We try to choose by induction on $\zeta < \lambda^+$
where $\lambda = |T| + |B|$ a sequence $\bar b^\zeta = \langle b_t:t
\in J_\zeta \rangle$ and $\eta_\zeta,\bar s_\zeta,\bar
t_\zeta,J^*,\varphi_\zeta,
\bar c_\varepsilon,\bar d_\varepsilon$
\mr
\item "{$(a)$}"   $J_\zeta$ is a linear order, increasing continuous
with $\zeta$
\sn
\item "{$(b)$}"  $J_0 = I_0,J_{\varepsilon +1} \backslash
J_\varepsilon$ is finite
\sn
\item "{$(c)$}"  $\bar b^\zeta$ is an indiscernible sequence over $A$
\sn
\item "{$(d)$}"  if $\zeta = \varepsilon +1$ then $n_\varepsilon <
\omega,\bar s_\varepsilon \in {}^{n_\varepsilon}(J_\zeta),\bar
t_\varepsilon \in {}^{n_\varepsilon}(J_\zeta),\varphi_\varepsilon =
\varphi_\varepsilon(\bar x_0,\dotsc,\bar x_{n_\varepsilon},\bar
c_\varepsilon,\bar d_\varepsilon),\bar c_\varepsilon \subseteq B,\bar
d_\varepsilon \subseteq A$ and $J^*_\varepsilon = \cup\{\bar
s_\xi \char 94 \bar t_\xi:\xi < \varepsilon\}$
\sn
\item "{$(e)$}"  $\bar s_\varepsilon \sim_{J^*_\varepsilon} \bar
t_\varepsilon$ and $\models \varphi[\bar b_{\bar s_\varepsilon},\bar
c_\varepsilon,\bar d_\varepsilon] \and \neg \varphi[\bar b_{\bar
t_\varepsilon},\bar c_\varepsilon,\bar d_\varepsilon]$ where $\bar
b_{\langle t_\ell:\ell < n \rangle} = \bar b_{t_0} \char 94 \bar
b_{t_1} \char 94 \ldots \char 94 \bar b_{t_{n-1}}$.
\ermn
If we succeed, \wilog \, $n_\varepsilon = n_*,\varphi_\varepsilon =
\varphi_*,
bar c_\varepsilon = \bar c^*$, and we get contradiction to
\scite{3.4}.  If we are stuck at stage $\varepsilon$, then $\bar{\bold
b}^\varepsilon$ is as required.
\enddemo
\bigskip

\remark{Concluding Remark}  We can define when an endless
indiscernible sequence is orthogonal to a set and the dimensional
independence property and prove natural properties, we intend to
pursue this.
\endremark

\newpage
    
REFERENCES.  
\bibliographystyle{lit-plain}
\bibliography{lista,listb,listx,listf,liste}

\shlhetal
\enddocument

\bye